\documentclass[11pt]{amsart}
\usepackage{latexsym,amssymb,amsmath,amscd,amsthm}
\numberwithin{equation}{section}
\theoremstyle{remark}

\newtheorem{proposition}{{\bf PROPOSITION}}[section]
\newcommand{\bq}{\begin{equation}}
\newcommand{\bea}{\begin{array}}
\newcommand{\eea}{\end{array}}

\newcommand{\ga}{\alpha}
\newcommand{\gep}{\epsilon}
\newcommand{\gD}{\Delta}
\newcommand{\gl}{\lambda}
\newcommand{\gL}{\Lambda}
\newcommand{\gb}{\beta}

\newcommand{\mf}{\mathfrak}
\newcommand{\mc}{\mathcal}

\newcommand{\wg}{\wedge}

\newcommand{\ci}{\circ}
\newcommand{\ul}[1]{\underline{#1}}

\newcommand{\ola}{\overleftarrow}
\newcommand{\go}{\omega}
\newcommand{\gO}{\Omega}

\newcommand{\gt}{\theta}
\newcommand{\gs}{\sigma}
\newcommand{\gz}{\zeta}
\newcommand{\gag}{\gamma}
\newcommand{\gd}{\delta}
\newcommand{\pp}{\partial}

\newcommand{\ora}{\overrightarrow}

\newcommand{\tl}{\tilde}
\newcommand{\na}{\nabla}
\newcommand{\gk}{\kappa}

\newcommand{\bs}{\blacksquare}

\newcommand{\bgs}{\bigstar}

\title{DISCRETIZATION, MOYAL, AND INTEGRABILITY}
\author{Robert Carroll\\University of Illinois, Urbana, IL 61801}

\dedicatory{{\bf \`A la m\'emoire de Joan}}

\date{May, 2001 - email: rcarroll@math.uiuc.edu}

\begin{document}

\bibliographystyle{plain}

\begin{abstract} 
Connections of KP, qKP, and Moyal type dKP constructions are developed.  Some expansion
of the Moyal KP constructions of Kemmoku-Saito is given with clarification of the role
of spectral variables as a phase space.
\end{abstract}

\maketitle




\section{INTRODUCTION}
\renewcommand{\theequation}{1.\arabic{equation}}
\setcounter{equation}{0}

This is a kind of sequel to \cite{czj} with connections to integrable
systems.  It is partially expository, clarifying the exposition in
\cite{ch} based on work of Kemmoku and Saito \cite{kaw,kas,kbe,kzq},
and some new aspects are indicated concerning q-integrable systems
and Moyal deformations.  In \cite{czj} we showed how a certain phase
space discretization in \cite{kaw,kas,kbe} is related to a
q-discretization leading to a new q-Moyal type bracket.  We discuss
here, following \cite{kaw,kbe}, also an analogous q-discretization
for KP related functions $A(z,\gz)=\sum a_{mn}z^m\gz^n$ with spectral
variables $(z,\gz)$ as phase space entities.  The version in
\cite{ch} is expanded and clarified.  We also indicate some features
of q-KP following \cite{au,at,dan,hy,ia,tu} and describe connections
to dispersionless theories and Moyal following \cite{gal,sh}.

\section{BACKGROUND AND ORIENTATION}
\renewcommand{\theequation}{2.\arabic{equation}}
\setcounter{equation}{0}

We run in parallel here the phase space discretization of
\cite{kaw,kas,kbe}, as expanded in \cite{ch,czj}, together with the
KP variants of \cite{kaw,kbe} as described in \cite{ch} (with some
correction of typos, etc.).  Thus for the phase space picture
involving
${\bf x}=(x,p)$ one has ($f=f(x,p)$)
\bq\label{1t}
X_f^D=\int da_1da_2v_{\gl}[f](x,p,a_1,a_2)\nabla_{{\bf a}};\,\,\nabla_{{\bf a}}=\frac
{1}{\gl}Sinh(\gl\sum a_i\pp_i);
\end{equation}
$$v_{\gl}[f]=\left(\frac{\gl}{2\pi}\right)^2\int db_1db_2e^{-i\gl(a_1b_2-a_2b_1)}
f(x+\gl b_1,p+\gl b_2)$$
\bq\label{2t}
X_f^D=\frac{\gl^3}{4\pi^2}\int d{\bf a}\int d{\bf b}e^{-i\gl({\bf a}\times{\bf b})}
\nabla_{{\bf b}}f\nabla_{{\bf a}};\,\,
X_f^Dg=\{f,g\}_M;\,\,[X_f^D,X_g^D]h=X^D_{\{f,g\}_M}h
\end{equation}
A variation of this in a q-lattice was developed in \cite{czj} (cf. also \cite{dam})
involving
\bq\label{58}
X_f^D=\sum_{m,n}v_q[f](x,p,m,n)\check{\na}_{mn};\,\,v_q[f]=\sum_{r,s}
q^{ms-nr}f(q^rx,q^sp);
\end{equation}
$$\check{\na}_{mn}g=\frac{g(xq^m,pq^n)-g(xq^{-m},pq^{-n})}
{(q^m-q^{-m})(q^n-q^{n})xp};$$
$$X_f^Dg=\frac{1}{2xp}\sum_{m,n,r,s}q^{ms-nr}\frac{[f(q^rx,q^sp)-f(q^{-r}x,q^{-s}p)]
[g(q^mx,q^np)-g(q^{-m}x,q^{-n}p)]}{(q^m-q^{-m})(q^n-q^{-n})}$$
Here the lattice structure was given a priori (no connection to $\gl$) and the role of 
$\gl$ is played by $q^m-q^{-m},\,\,q^n-q^{-n},$ etc.
For the KP situation with 
${\bf (P1)}\,\,A(z,\gz)=\sum a_{mn}z^m\gz^n\,\,(m,n\in{\bf Z})$ one writes
\bq\label{3t}
X_A^D=\int d{\bf a}v_{\gl}[A]\nabla_{{\bf a}};\,\,\nabla_{{\bf
a}}\sim\frac{1}{\gl}Sin(\gl{\bf a}\cdot\pp);\,\,
v_{\gl}[A]=\frac{\gl^2}{4\pi^2}\int d{\bf b}e^{-i\gl({\bf a}\times {\bf b})}
A({\bf x}+i\gl{\bf b})
\end{equation}
where $z\sim e^p$ and $\gz\sim e^x$ with ${\bf x}\sim (x,p)$ so for $A_{mn}=z^m\gz^n$
one has ${\bf (P2)}\,\,A_{mn}({\bf x}+i\gl{\bf b})=e^{mp}e^{nx}exp(i\gl(mb_2+nb_1)=exp
[m(p+i\gl b_2)]exp[n(x+i\gl b_1)]$ (note $\pp_x=\gz\pp_{\gz}=\pp_{log(\gz)}$, etc.). 
Then (specifying $\nabla_{mn}\sim\nabla_{m,-n}$)
\bq\label{4t}
v_{\gl}[A_{mn}]=A_{mn}\gd(m-a_1)\gd(n+a_2);
\end{equation}
$$X^D_{A_{mn}}=A_{mn}\nabla_{mn}=A_{mn}
\frac{1}{\gl}Sin[\gl(m\pp_{log(\gz)}-n\pp_{log(z)}]$$
Further, for ${\bf (P3)}\,\,X_{mn}^D=z^m\gz^n\nabla_{mn}\sim z^m\gz^n\nabla_{m,-n}$,
one has
\bq\label{5t}
X_A^D=\sum a_{mn}z^m\gz^n\nabla_{mn}=\sum A_{mn}\nabla_{mn}=\sum a_{mn}X^D_{A_{mn}};
\end{equation}
$$[X_{mn}^D,X_{pq}^D]=\frac{1}{\gl}Sin[\gl(np-mq)]X^D_{m+p,n+q}$$
and for $f,g\sim A,B$ one has 
${\bf (P4)}\,\,[X_A^D,X_B^D]=-X_{\{A,B\}_M}$ with
\bq\label{6t}
\{A,B\}_M=\frac{-1}{\gl}Sin\left[\gl(\pp_{log(z_1)}
\pp_{\log(\gz_2)}-\pp_{log(\gz_1)}\pp_{log(z_2)})\right]
A(z_1,\gz_1)B(z_2,\gz_2)|_{(z,\gz)}
\end{equation}
\\[3mm]\indent
Next recall from \cite{ch,kaw}
\bq\label{7t}
P_{F_w}=\frac{\hbar}{4}\int d{\bf a}\int d{\bf b}e^{i\hbar({\bf a}\times {\bf
b})}F_w({\bf x}+(\hbar/2){\bf b})\gD^{{\bf a}}
\end{equation}
where $<\gD^{{\bf b}},\nabla_{{\bf a}}>=\gd_{ba}$ and $F_w$ is the Wigner function
\bq\label{8t}
F_w(p,x)=\frac{1}{2\pi}\int dy\phi\left(x+\frac{\hbar}{2}y\right)\phi^*\left(x-\frac
{\hbar}{2}y\right)e^{-ipy}
\end{equation}
Then for $X_f^Q=\hbar X_f^D$ one has ${\bf (P5)}\,\,<P_{F_w},X^Q_f>=\int dxdp
F_w(x,p)f(x,p)=<\hat{f}>$.  For KP we have from \cite{ch,kaw}
\bq\label{9t}
F_{KP}(z,\gz)=\int dx\sum_{\ell\in{\bf
Z}}\psi(q^{\ell/2}z)\psi^*(q^{-\ell/2}z)\gz^{-\ell}=\sum f_{mn}z^{-m}\gz^{-n}
\end{equation}
Let now ${\bf (P6)}\,\,<\gD^{mn},\nabla_{pq}>=\gd_{mp}\gd_{nq}$ and write for 
$A=\sum a_{mn}z^m\gz^n$
\bq\label{10t}
\gO_{F_{KP}}(z,\gz)=\sum_{mn}f_{mn}z^{-m}\gz^{-n}\gD^{mn};
\end{equation}
$$<\gO_{F_{KP}},X_A^D>=-\oint\frac{dz}{2\pi iz}\oint\frac{d\gz}{2\pi i\gz}
F_{KP}(z,\gz)A(z,\gz)=\sum f_{mn}a_{mn}$$
In particular for $X_{mn}^D=z^m\gz^n\nabla_{mn}$ one has
\bq\label{11t}
f_{mn}=<\gO_{F_{KP}},X^D_{mn}>=-\int dx\oint\frac{dz}{2\pi iz}z^m\psi(q^{n/2}z)
\psi^*(q^{-n/2}z)
\end{equation}
\\[3mm]\indent
We now recall the Orlov-Schulman operators (cf. \cite{cqq,dap,daq,ot}),
i.e. ${\bf (P7)}\,\,z^m\pp_z^{\ell}\phi=M^{\ell}L^m\phi$ with flow equations
$\pp_{m\ell}\phi=-(M^{\ell}L^m)_{-}\phi$ for $\phi=Wexp(\xi)$ (i.e.
$z^m\pp_z^{\ell}\sim\pp_{m\ell}$).  Recall also $L\phi=z\phi$ and
$\pp_n\phi=L^n_{+}\phi$; from this ${\bf (P8)}\,\,(\pp_n-L^n_{+})W=W(\pp_n-\pp^n)$ and
$\pp_nW=(\pp_nW)+W\pp_n$ implies $(\pp_nW)W^{-1}=L^n_{+}-W\pp^nW^{-1}=-L^n_{-}$ or
$\pp_nW=-L^n_{-}W$ (Sato equation).  Then one can rewrite the right side $\Xi$ of
\eqref{11t} in terms of the $\pp_{m\ell}$. Indeed $q^{nz\pp_z}\psi(z)=\psi(q^nz)$ and
under a change of variables $z\to q^{n/2}z$ the integrand in \eqref{11t} becomes
$(dz/2\pi iz)(q^{n/2}z)^m\psi(q^nz)\psi^*(z)$ so
\bq\label{12t}
\Xi=-\int dx\oint\frac{dz}{2\pi iz}\left(z^mq^{n[z\pp_z+(m/2)]}\psi(z)\right)\psi^*(z)
\end{equation}
One can write then for $\gl=log(q)$ a formula ${\bf (P9)}\,\,q^{nz\pp_z}=exp(n\gl
z\pp_z)=\sum_0^{\infty}(n\gl)^j(z\pp_z)^j/j!$ so
\bq\label{13t}
\Xi=-\int dx q^{mn/2}\oint\frac{dz}{2\pi
iz}z^m\left(\sum_0^{\infty}\frac{(n\gl)^j}{j!}z^j\pp_z^j\psi(z)\right)\psi^*(z)=
\end{equation}
$$=-\int dx q^{nm/2}\sum_{j=0}^{\infty}\frac{(n\gl)^j}{j!}\oint\frac{dz}{2\pi i}\left(
z^{m+j-1}\pp_z^j\psi(z)\right)\psi^*(z)$$
Now from {\bf (P7)} we can write ${\bf (P10)}\,\,\pp_{m\ell}\psi=-(M^{\ell}L^m)_{-}\psi$
with $z^m\pp_z^{\ell}\psi=M^{\ell}L^m\psi$ while in \cite{kaw} we have
${\bf (P11)}\,\,\pp_{m+\ell,\ell}\pp log(\tau)=\oint[dz/2\pi i](z^{m+\ell}\pp_z^{\ell}
\psi(z))\psi^*(z)= Res(z^{m+\ell}\pp_z^{\ell}\psi(z))\psi^*(z)$.  There are various 
formulas in this direction and we recall some results from \cite{cqq,dap,daz,day,mss,
ow}.  Thus one has a lemma of Dickey that for $P=\sum p_k\pp^{-k}$ and $Q=\sum
q_k\pp^{-k}$ follows $Res_{\gl}(Pexp(\gl x)Qexp(-\gl x)=Res_{\pp}PQ^*$ (here
$x^*=x,\,\,\pp^*=-\pp$, and $(AB)^*=B^*A^*$).  Now set
$W=1+\sum_1^{\infty}w_j\pp^{-j}$ and $(W^*)^{-1}=1+\sum_1^{\infty}w_j^*\pp^{-j}$.
The flow equations are $\pp_mu=\pp Res(L^m)=\pp Res(L^m_{-})$ corresponding,
for $w=w_1$ (with $u=-\pp w$), to $-\pp_m\pp w=\pp Res(L^m)$ or $\pp_mw=-Res(L^m)=
Res(L^m_{-})$.  One knows also ${\bf
(P12)}\,\,Res(M^nL^{m+1})=Res(z^{m+1}\pp_z^n\psi)\psi^*$.  Further
$\psi\psi^*=1-\sum_1^{\infty}[\pp_iw/z^{i+1}]$ with $\pp_1\sim\pp_x,\,\,\pp_2\sim\pp_y,
\,\,\pp_3\sim\pp_t,\cdots$.  This corresponds to $\psi\psi^*=\sum s_nz^{-n}$ with 
$\pp s_{n+1}=\pp_nu\sim\pp s_{n+1}=-\pp_n\pp w\sim s_{n+1}=-\pp_nw$.  Thus one
writes
\bq\label{14t}
Res(z^m\pp_z^n\psi)\psi^*=Res M^nL^m
\end{equation}
and in \cite{kaw} $\pp_{k\ell}$ is defined via ${\bf (P13)}\,\,\pp_{k\ell}W=
-(M^{\ell}L^k)_{-}W$ (additional Sato equations - cf. {\bf (P8)}).  Consequently
using the Dickey lemma for $w=w_1$ and {\bf (P13)} with $L\psi=z\psi$ and
$\pp_z\psi=M\psi$ ($[L,M]=1$) one looks at $\pp_{k\ell}w=-Res_{\pp}(M^{\ell}L^k)$.  Set 
$P=Wx^{\ell}\pp^k$ and $Q=(W^*)^{-1}$ so $Q^*=W^{-1}$ and
$PQ^*=Wx^{\ell}\pp^kW^{-1}=Wx^{\ell}W^{-1}W\pp^kW^{-1}$.  Now one can define $L,M$ also
via $L=W\pp W^{-1}$ and $M=W(\sum kt_k\pp^{k-1})W^{-1}$ so for $\psi=wexp(\xi)$ we have
$M\psi=W(\sum kt_k\pp^{k-1})exp(\xi)=W(\sum
kt_kz^{k-1})exp(\xi)=W\pp_zexp(\xi)=\pp_z\psi$ and $L\psi=W\pp exp(\xi)=z\psi$.
Hence
($e^{\pm\xi}$ replaces $e^{\pm zx}$ in an obvious manner where $\xi=\sum t_nz^n$)
\bq\label{15t}
\pp_{k\ell}w=-Res_{\pp}M^{\ell}L^k=-Res_z(Pe^{\xi})(Qe^{-\xi})=
-Res_z(Wx^{\ell}\pp^ke^{\xi}
(W^*)^{-1}e^{-\xi})=
\end{equation}
$$=-Res_z[z^kW(x^{\ell}e^{\xi})(W^*)^{-1}e^{-\xi}]=
-Res_z\left(z^k[\pp_z^{\ell}\psi(z)]\psi^*(z)\right)$$
Now in \cite{ch} there is some confusion about using $\oint dz/2\pi i$ or $\oint
dz/2\pi iz$ as a symbol for a residue.  If the former is adopted consistently then
\eqref{13t} becomes (via \eqref{15t})
\bq\label{16t}
\Xi=-\int dx q^{mn/2}\sum_{j=0}^{\infty}\frac{(n\gl)^j}{j!}Res\left[\left(z^{m+j-1}
\pp_z^j\psi(z)\right)\psi^*(z)\right]=
\end{equation}
$$=\int dx q^{mn/2}\sum_{j=0}^{\infty}\frac{(n\gl)^j}{j!}\pp_{m+j-1,j}w$$
where $u=-\pp w\Rightarrow w=-\pp^{-1}\pp^2log(\tau)=-\pp log(\tau)$ which leads to
\bq\label{17t}
\Xi=\int dx {\mc D}_{mn}\left(-\pp log(\tau)\right);\,\,{\mc D}_{mn}=q^{mn/2}\sum_
{j=0}^{\infty}\frac{(n\gl)^j}{j!}\pp_{m+j-1,j};
\end{equation}
$$f_{mn}=<\gO_{F_{KP}},X_{mn}^D>=\int dx{\mc D}_{mn}(-\pp log(\tau))$$
Thus actually it is $F_{KP}$ in \eqref{9t} which generates ${\mc D}_{mn}$ via $X_{mn}^D$
(which eliminates the $a_{mn}$ in \eqref{9t}).  One can write now 
${\mc D}_A=\sum a_{mn}{\mc D}_{mn}$ with 
\bq\label{18t}
<\gO_{F_{KP}},X^D_A>=\int dx {\mc D}_A(-\pp log(\tau))=\tl{A}(t);
\end{equation}
$$\pp_r\tl{A}=-\int dx{\mc D}_A(\pp_r\pp log(\tau))=-\int dx{\mc D}_AJ_r$$
where $J_r=\pp_r\pp log(\tau)$ is a first integral of KP (i.e. $H_r=\int dx J_r$ is a
Hamiltonian).  Here one recalls that $s_{n+1}=\pp^{-1}\pp_nu=\pp_n\pp og(\tau)=J_n$ are
conserved densities (cf. \cite{cqq} and note that the $\pp s_{n+1}=K_{n+1}=\pp_nu$ are
symmetries determining the standard KP flows).  Note also that the ${\mc D}_{mn}$ flows
are independent and commuting with the $\pp_r$ so ${\mc D}_A$ is independent and
commutes with $\pp_r$.  Now one defines ${\bf (P14)}\,\,-\int dx
{\mc D}_AJ_r={\mc D}_{\tl{A}}\cdot H_r={\mc D}_{\tl{A}}\cdot\int dx J_r$ where 
${\mc D}_{\tl{A}}$ corresponds to a Hamiltonian vector field (e.g. $X_f$ acting as
$X_fg$ in \eqref{2t}); i.e. one writes
\bq\label{19t}
\pp_r\tl{A}=\{\tl{A},H_r\}_M^{KP}
\end{equation}
This is rather a stretch of imagination but perhaps morally correct at least since it is
consistent with the Heisenberg notation ${\bf (P15)}\,\,dA/dt=\{A,H\}_M$.  Here
$\{\tl{A},H_r\}_M^{KP}$ makes no recourse to phase space variables $z,\gz$ however and
any relation to e.g. \eqref{6t} is vague at best.  As in \cite{ch} we can say however
that $\pp_r\tl{A}={\mc D}_{\tl{A}}\cdot H_r$ has the structure $\pp_r<\gO_{F_{KP}},
X_A^D>=\pp_r\tl{A}$ and one imagines e.g. a Heisenberg picture
\bq\label{20t}
\pp_tX_A^D=-[X_A^D,X_H^D]=X_{\{A,H\}_M}\equiv \pp_t<\gO_{F_w},X_A^D(t)>=<\gO_{F_w},
X^D_{\{A(t),H\}_M}>
\end{equation}
(cf. also \cite{ch,czj}).
\\[3mm]\indent
{\bf REMARK 2.1.}  The q-lattice version \eqref{58} of phase space discretization was
developed following an analogue of Fourier transform techniques dealing with 
\eqref{1t} - \eqref{2t}.  On the other hand \eqref{4t}, derived via $z\sim e^p$ and
$\gz\sim e^x$, leads for $q\sim exp(i\gl)$ to
\bq\label{36t}
X^D_{A_{mn}}B(\gz,z)=\frac{1}{2i\gl}[B(q^m\gz,q^{-n}z)-B(q^{-m}\gz,q^nz)]
\end{equation}
which has some similarity to $\check{\nabla}_{mn}g$ of \eqref{58}.
In accord with procedures in \cite{czj} one might expect here a directive to modify
$\gl$ in the denominator of \eqref{36t} via e.g. $(q^m-q^{-m})(q^n-q^{-n})
z\gz$.  At first
sight one is tempted to look for $\gz,z$ as arising in vertex operators (cf.
\cite{au,al,cw,cqq,cn,dap,mzz})
\bq\label{37t}
X(z,\gz,t)\tau=e^{\sum_1^{\infty}t_k(\gz^k-z^k)}e^{\sum_1^{\infty}(z^{-k}-\gz^{-k})
\pp_k/k}\tau=
\sum_{m=0}^{\infty}\frac{(\gz-z)^m}{m!}\sum_{p=-\infty}^{\infty}z^{-p-n}W_p^m\tau
\end{equation}
(where the $W_p^m$ are expressed in terms of currents 
$J_k^v$ involving time derivatives $\pp_n$) leading perhaps to a new perspective
on vertex operators. 
The dynamical analogies are not immediately clear
but an origin of spectral variables as phase space variables is suggested in Remark 4.2
via the anti-isomorphism of PSDO and the z-operators of Dickey (cf. \cite{daz,day}).
Nevertheless we give here for completeness a few more formulas involving vertex
operators. 
Thus one can imagine of course a
Taylor type expansion (cf. \cite{al})
\bq\label{38t}
X(z,\gz,t)=\sum_{m=0}^{\infty}\frac{(\gz-z)^m}{m!}\pp_{\gz}^mX(z,\gz,t)|_{\gz=z};\,\,
\pp_{\gz}^mX(z,\gz,t)|_{\gz=z}=\sum_{p=-\infty}^{\infty}z^{-p-m}W_p^m
\end{equation}
which would suggest that for $n\in{\bf Z}$
\bq\label{39t}
Res[z^n\pp_{\gz}^mX(z,\gz,t)]|_{\gz=z}=W^m_{n-m+1}
\end{equation}
(i.e. $n-p-m=-1$) and for $n=m+s$ this gives $W_{s+1}^m$.  A variation on this
(Adler-Shiota-van Moerbeke theorem) involves ${\bf
(P16)}\,\,Res[z^{m+s}\pp_z^mX(z,\gz,t)|_{\gz=z}\tau\sim W_s^{m+1}\tau/(m+1)$ or more
generally (cf. \cite{au,dap,mzz})
\bq\label{40t}
-\frac{(M^mL^{m+s})_{-}\psi}{\psi}\sim
\frac{Y_{m+s,m}\psi}{\psi}=(e^{\eta}-1)\frac{W_s^{m+1}\tau}{(m+1)\tau}
\end{equation}
where $[exp(-\eta)-1]f(t)=f(t-[z^{-1}])-f(t)$ via $\eta=\sum_1^{\infty}\pp_j/jz^j$
(note in \eqref{40t} the last term requires some $z$ dependence - further in taking
residues only a $D(W_s^{m+1}\tau/\tau)$ term arises).
A possible way now to envision $z,\gz$ as phase space variables might be in
terms of 
$\gz\sim\pp_z$ or $z\sim L$ and $\gz\sim M$.  We recall also from \cite{al} that the
generating functions
\bq\label{41t}
W_{\gl}^v=\sum_{-\infty}^{\infty}\gl^{-p-v}W_p^v;\,\,J_{\gl}^v=\sum_{-\infty}^{\infty}
\gl^{-p-v}J_p^v
\end{equation}
can be considered as stress energy tensors and ($D\sim\pp_x$)
\bq\label{42t}
\psi^*(\gl,t)\psi(\mu,t)=\sum_1^{\infty}\frac{(\mu-\gl)^{j-1}}{j!}D\frac{W_{\gl}^j
(\tau)}{\tau}=\frac{1}{\mu-\gl}D\left(\frac{X(\gl,\mu,t)\tau}{\tau}\right)
\end{equation}
\bq\label{43t}
v\psi^*(\gl,t)\pp_{\gl}^{v-1}\psi(\gl,t)=D\left(\frac{1}{\tau}\sum_{-\infty}^{\infty}
\gl^{-p-v}W_p^v(\tau)\right)
\end{equation}
Hence from \eqref{43t}
\bq\label{44t}
Res[\gl^{\ell+v}\psi^*(\gl,t)\pp_{\gl}^v\psi(\gl,t)]=
\end{equation}
$$=Res\left[D\frac{1}{(v+1)\tau}\left(
\sum_{-\infty}^{\infty}\gl^{-p+\ell-1}W_p^{v+1}(\tau)
\right)\right]=D\left(\frac{W_{\ell}^{v+1}(\tau)}
{(v+1)\tau}\right)$$
which is consistent with \eqref{40t}.
$\hfill\bs$
\\[3mm]\indent
In connection with KP and Moyal we add a few results from \cite{gal}
(cf. also \cite{ch,sh,tz}).
Thus let $\tilde{M}=M\oplus T$ where $T\sim$ times $\{t_n\}$ and for some
$*$ product on $M$ write
\bq
u(x,t)* v(x,y)=
\left.exp\left(\kappa\omega^{ij}\frac{\partial}{\partial x^i}\frac
{\partial}{\partial \tilde{x}^j}\right)u(x,t)v(\tilde{x},t)\right|_{\tilde{x}=x}
\label{X12}
\end{equation}
(note $X_f=\omega^{ij}(\partial f(x,t)/\partial x^i)(\partial/\partial x^j)$
will be time dependent).  Let $L=\partial+\sum_1^{\infty}u_n(x,t)
\partial^{-n}$ be the Lax operator for KP. 
One then applies the geometrical framework to obtain a Moyal KP
hierarchy $KP_{\kappa}$, based on deformation of dKP, which is 
equivalent to the Sato hierarchy based on PSDO.
Similar calculations apply to Toda and dToda, KdV and dKdV, etc.
Further the geometrical picture can be phrased in the $Sdiff_2$ format
with ${\mc L},\,{\mc M}$ etc.  It seems from this that if one starts
with dKP as a basic Hamiltonian system with Hamiltonians
${\mc B}_n$ and standard P brackets then $KP_{\kappa}$ can be considered
as a quantization of dKP is some sense with quantum integrals of motion
$B_n(\kappa)$ which for $\kappa=1$ say is equivalent to KP (or $\kappa=
1/2$ in \cite{gal}).  The $B_n(\kappa)$ would perhaps have to be extracted
from KP after establishing the isomorphism (cf. \cite{gal}) and we turn
briefly to this approach now.  Thus in \cite{gal} one writes the Sato KP
hierarchy via ($v_{-2}=1,\,v_{-1}=0$)
\bq
\partial_mL=[L^m_{+},L]\,\,(m\geq 1);\,\,L=\sum_{-2}^{\infty}v_n
(\tilde{x})\partial_x^{-n-1}
\label{X15}
\end{equation}
for $\tilde{x}=(x,t_2,\cdots$) while the Moyal KP hierarchy is written via
($u_{-2}=1,\,u_{-1}=0$)
\bq
\Lambda=\sum_{-2}^{\infty}u_n(\tilde{x})\lambda^{-n-1};\,\,\partial_m
\Lambda=\{\Lambda^m_{+},\Lambda\}_M\,\,(m\geq 1)
\label{X16}
\end{equation}
where $\Lambda^m_{+}\sim (\Lambda^{*m})_{+}$ with
\bq
f*g=\sum_0^{\infty}\frac{\kappa^s}{s!}\sum_{j=0}^s(-1)^j{s \choose j}
(\partial_x^j\partial_{\lambda}^{s-j}f)(\partial_x^{s-j}\partial_{\lambda}^jg)
\label{X17}
\end{equation}
leading to
\bq
\{f,g\}_{\kappa}=\sum_0^{\infty}\frac{\kappa^s}{(2s+1)!}\sum_{j=0}^{2s+1}
(-1)^j{2s+1 \choose j}(\partial_x^j\partial_{\lambda}^{2s+1-j}f)
(\partial_x^{2s+1-j}\partial_{\lambda}^jg)
\label{X18}
\end{equation}
Note $lim_{\kappa\to 0}\{f,g\}_{\kappa}=\{f,g\}=f_{\lambda}g_x-
f_xg_{\lambda}$ so $(KP)_M\to dKP$ as $\kappa\to 0$, namely
$\partial_m\Lambda=\{\Lambda^m_{+},\Lambda\}$ with $\Lambda^m\sim
\Lambda\cdots\Lambda$.  The isomorphism between $(KP)_{Sato}$ and $(KP)_M$ 
is then determined by relating $v_n$ and $u_n$ in the form $(\kappa=1/2$)
\bq
u_n=\sum_0^n2^{-j}{n \choose j}v^j_{n-j}
\label{X19}
\end{equation}
where $n=0,1,\cdots$ and $v^j=\partial_x^jv_0$.
\\[3mm]\indent
Let us recall the $KP_{\gk}$ theory of \cite{gal,sh} where dKP is built up as
follows (cf. \cite{ch,cn,co,tz}). 
We will assume some familiarity with KP and dKP as in \cite{ch} and only
recall formulas (extensive references are given in \cite{ch}).  Then the KP Lax
operator has the form
${\bf (P17)}\,\,L=\pp+\sum_1^{\infty} u_{n+1}\pp^{-n}$ with $u_2=u$ and
$u_i=u_i(x,x_n)$ where $x_1\sim x$.
The $x_n$ (or equivalently $t_n$) for $n\geq 2$ correspond to time variables with flows
${\bf (P18)}\,\,\pp_nL=[B_n,L]$ for $B_n=L_{+}^n$.  There is a dressing or gauge
operator ${\bf (P19)}\,\,W=1+\sum_1^{\infty}w_n\pp^{-n}$ determined via $L=W\pp
W^{-1}\,\,(\pp=\pp_x)$.  For wave functions or Baker-Akhiezer (BA) functions 
$\psi(x,\gl)=Wexp(\xi),\,\,\xi=\sum_1^{\infty}x_n\gl^n$ one hs $L\psi=\gl\psi$ and
$\pp_m\psi=B_m\psi$.  Further $\psi^*=(W^*)^{-1}exp(-\xi)$ and $L^*\psi^*=\gl\psi^*$
with $\pp_m\psi^*=-B_n^*\psi^*$ (here $B_n^*=(L^*)^n_{+}$ and $\pp^*=-\pp$).  The
equation ${\bf (P20)}\,\,\pp_nW=-L_{-}^nW$ is called the Sato equation.  Concerning
$\gl$ derivatives one has
\bq\label{2p}
\psi_{\gl}=W\left(\sum_1^{\infty}kx_k\gl^{k-1}\right)exp(\xi)=M\psi;\,\,M=W
\left(\sum_1^{\infty}kx_k\pp^{k-1}\right)W^{-1}
\end{equation}
and $[L,M]=1$.  The tau function arises in a vertex operator equation (VOE)
\bq\label{3p}
\psi(x,\gl)=\frac{X(\gl)\tau}{\tau}=e^{\xi}\frac{\tau_{-}}{\tau}=\frac{e^{\xi}}{\tau}
\tau\left(x_j-\frac{1}{j\gl^j}\right)
\end{equation}
The Hirota bilinear identity is ${\bf
(P21)}\,\,0=\oint_C\psi(x,\gl)\psi^*(y,\gl)d\gl$, where $C$ is a circle around
$\gl=\infty$, and this leads to various Hirota bilinear formulas.  In particular one
has $(\tl{\pp}\sim (\pp_j/j))$
\bq\label{4p}
\sum_0^{\infty}p_n(-2y)p_{n+1}(\tl{\pp})exp\left(\sum_1^{\infty}y_j\pp_j\right)\tau\cdot
\tau=0
\end{equation}
where the $p_j$ are Schur polynomials and ${\bf (P22)}\,\,\pp_j^ma\cdot b=(\pp^m/\pp
s_j^m)a(t_j+s_j)b(t-g-s_j)|_{s_j=0}$.
The KP equation is included in \eqref{4p} in the form ${\bf
(P23)}\,\,(\pp^4+3\pp_2^2-4\pp_1\pp_3)\tau\cdot\tau=0$. 
\\[3mm]\indent
{\bf REMARK 2.2.}   
For dKP traditionally one
thinks of fast and slow variables $\gep t_i=T_i$ (shifting now $x_i\to t_i$ for $i\geq
2$ with $x_1\sim x$) and $\pp_n\to\gep\pp/\pp T_n$ with $\pp=\pp_x\to\gep\pp_X$ and
$\pp_x^{-1}\to (1/\gep)\pp_X^{-1}$.  Then one writes $u_i(x,t_n)\to\tl{u}_i(X,T_n)$
and this passage (where one usually assumed $u_i(X/\gep,T_n/\gep)=
\tl{u}(X,T_n)+O(\gep)$ has always seemed unrealistic; however in certain situations it
is perfectly reasonable (see e.g. \cite{ch,cm} under $(X,\psi)$ duality).  A priori
if one simply substitutes in a power series ${\bf (P24)}\,\,\sum
a_{\ga}x^{\ga_1}t_2^{\ga_2}\cdots\to \sum a_{\ga}\gep^{-\sum
\ga_i}X^{\ga_1}T_2^{\ga_2}\cdots$ there will be horrible divergences as $\gep\to 0$ so
one is led to think of sums of simple homogeneous functions (e.g. $f(x,t_n)\to
f(X/\gep,T_n/\gep)=\sum_0^{\infty}\gep^jF_j(X,T_n)$ with sums as in {\bf (P24)}
yielding terms with $-\sum\ga_i\geq 0$.  Note that this can be achieved for arbitrary
powers of x and a finite number of $t_n\,\,(2\leq n\leq N)$ by insertion of some
suitably large negative power of say $t_{N+1}$ in each monomial; then one could worry
about the meaning of $t_{N+1}$ later.  Perhaps this is an argument for some ultimate
projectivization via $T_{N+1}$ corresponding to some universal time.$\hfill\bs$
\\[3mm]\indent
In any event for dKP one writes $(T\sim (X,T_n)$ - cf. \cite{ch,cn,co,tz})
\bq\label{5p}
L_{\gep}=\gep\pp+\sum_1^{\infty}u_{n+1}(\gep,T)(\gep\pp)^{-n};\,\,M_{\gep}=
\sum_1^{\infty}nT_nL_{\gep}^{n-1}+\sum_1^{\infty}v_{n+1}(\gep,T)L_{\gep}^{-n-1}
\end{equation}
with $u_{n+1}(\gep,T)=U_{n+1}(T)+O(\gep)$ and $v_{n+1}(\gep,T)=V_{n+1}(T)+O(\gep)$. 
Then set
\bq\label{6p}
\psi=\left[1+O\left(\frac{1}{\gl}\right)\right]e^{\sum_1^{\infty}(T_n/\gep)\gl^n}
=exp\left(\frac{1}{\gep}S(T,\gl)+O(1)\right);
\end{equation}
$$\tau=exp\left(\frac{1}{\gep^2}F(T)+O\left(\frac{1}{\gep}\right)\right)$$
Setting $P=S_X$ with ${\bf (P25)}\,\,\tau(T-(1/n\gl^n))exp[\sum_1^{\infty}T_n\gl^n]/
\tau(T)$ one obtains
\bq\label{7p}
\gl=P+\sum_1^{\infty}U_{n+1}P^{-n};\,\,P=\gl-\sum_1^{\infty}P_i\gl^{-i};\,\,
\mf{M}=\sum_1^{\infty}nT_n\gl^{n-1}+\sum_1^{\infty}V_{n+1}\gl^{-n-1}
\end{equation}
along with ${\bf (P26)}\,\,B_n=\sum_0^n b_{nm}\pp^m\to \mf{B}_n=\sum_0^n
b_{nm}P^m\sim\gl_{+}^n$ leading to (note $\pp_nM=[B_n,M]$)
\bq\label{8p}
\pp_n\gl=\{\mf{B}_n,\gl\};\,\,\pp_n\mf{M}=\{\mf{B}_n,\mf{M}\};\,\,\{\gl,\mf{M}\}=1
\end{equation}
where $\{A,B\}=\pp_PA\pp_XB-\pp_XA\pp_PB$.
\\[3mm]\indent
Now return to \cite{gal,sh} (sketched briefly in \eqref{X12} - \eqref{X19}).  Note
$\sum_{-2}^{\infty}v_n(\tl{x})\pp_x^{-n-1}=\pp+v_0\pp^{-1}+\cdots$ corresponds to
{\bf (P17)} with index changes and $v\sim u$, while
$\gL=\sum_{-2}^{\infty}u_n(\tl{x})\gl^{-n-1}$ corresponds to
$\gL=\gl+u_0\gl^{-1}+\cdots$ which corresponds to \eqref{7p} with $\gl\sim P$.  Thus
we see that the Moyal bracket in \eqref{X18} for example involves $X$ and $P$
derivatives where $(X,P)$ is the natural phase space for dKP.  The connection between
$(KP)_{Sato}$ and $(KP)_{\gk}$ for $\gl=1/2$ is obtained by relating $u_n\sim U_{n+1}$
in dKP to $v_n\sim u_{n+1}$ in KP.  It is important to realize here that no scaling is
involved in \cite{gal,sh}.  In any case
one can formulate the KP hierarchy as a quantization of dKP under the Moyal bracket.
The actual correspondence \eqref{X19} is not important here 
(see also below) and one could simply define
KP as $(KP)_{\gk}$ for $\gk=1/2$ and express it through phase space $(X,P)$ Moyal
brackets.  We will discuss below similar correspondences for q-KP and dKP under
suitable q-Moyal type brackets.  In this direction, following \cite{sh}, one would
have
\bq\label{9p}
\frac{\pp\gl}{\pp t_n}=\{\mf{B}_n,\gl\}_{\gk};\,\,\mf{B}_n\sim (\gl*\cdots
*\gl)_{+};\,\,\frac{\pp\mf{B}_n}{\pp t_m}-\frac{\pp\mf{B}_m}{\pp t_n}+\{\mf{B}_n,
\mf{B}_m\}_{\gk}=0
\end{equation}

\section{REMARKS ON Q-KP}
\renewcommand{\theequation}{3.\arabic{equation}}
\setcounter{equation}{0}

There are various approaches to q-KP and we mention e.g.
\cite{at,az,fam,fal,fa,faj,hy,ia,kah,kag,kaz,smj,tu,va,za}.
We will not dwell upon q-nKdV or q-NLS here, nor upon discrete KP as in
\cite{at,az,dan}.  Let us rather follow \cite{ia,tu} at first in writing
${\bf (P27)}\,\,Df(x)=f(xq)$ with $D_qf(x)=[f(xq)-
f(s)]/x(q-1)$ and we recall
\bq\label{9h}
(a;q)_0=1;\,\,(a;q)_k=\prod_0^{k-1}(1-aq^s);\,\,\left(\begin{array}{c}
n\\
k\end{array}\right)_q=\frac{(1-q^n)(1-q^{n-1})
\cdots(1-q^{n-k+1})}{(1-q)(1-q^2)\cdots(1-q^k)}
\end{equation}
Then using ${\bf
(P28)}\,\,D_q^n(fg)=\sum_0^n\binom{n}{k}_q(D^{n-k}D_q^kf)D_q^{n-k}g$
one obtains a formula ${\bf
(P29)}\,\,D_q^nf=\sum_0^{\infty}\binom{n}{k}_q(D^{n-k}D_q^kf) D_q^{n-k}\,\,(n\in{\bf
Z}$ and this is
shown to be correct via ${\bf (P30)}\,\,D_q^m(D_q^nf)= D_q^{m+n}f$ (proved in
\cite{ia}).  This leads to the formal adjoint to $P=
\sum a_iD_q^i$ as ${\bf (P31)}\,\,P^*=\sum(D_q^*)^ia_i$ where $D^*_q=-(1/q)D_{1/q}$
and the result that ${\bf (P32)}\,\,(PQ)^*=Q^*P^*$. Now consider the formal
q-PSDO ${\bf (P33)}\,\,L=D_q+a_0+\sum_1^{\infty}a_iD_q^{-i}$ leading to the
q-deformed KP hierarchy $(\pp L/\pp t_j)=[(L^j)_{+},L]$ (this differs by a factor of
$x(1-q)$ from the definitions in \cite{fa}).  Let S be the PSDO ${\bf(P34)}\,\,
S=1+\sum_1^{\infty}w_kD_q^{-k}$ satisfying $L=SD_qS^{-1}$ (whose existence is
proved easily).  The vector fields $\pp/\pp t_j$ can be extended via ${\bf
(P35)}\,\,(\pp S/\pp t_j)=-(L^j)_{-}S$ and will remain commutative. 
One uses now a nonstandard definition
\bq\label{10h}
exp_q(x)=\sum_0^{\infty}\frac{(1-q)^kx^k}{(q;q)_k}=exp\left(\sum_1^{\infty}
\frac{(1-q)^kx^k}{k(1-q^k)}\right)
\end{equation}
and $D_q^kexp_q(xz)=z^kexp_q(xz)$.  One uses also the notation ${\bf (P36)}\,\,
P|_{x/t}=\sum p_i(x/t)t^iD_q^i$ when $P=\sum p_iD^i_q$.  The q-wave function $w_q$
and its adjoint $w_q^*$ are defined via
\bq\label{11h}
w_q(x,t)=Sexp_q(xz)exp\left(\sum_1^{\infty}t_iz^i\right);\,\,w^*_q=(S^*)^{-1}|_
{x/q}exp_{1/q}(-xz)exp\left(-\sum_1^{\infty}t_kz^k\right)
\end{equation}
One can easily show (as in the classical case)
\bq\label{12h}
Lw_q=zw_q;\,\,\pp_mw_q=(L^m)_{+}w_q;\,\,L^*|_{x/q}w^*_q=zw_q^*;\,\,\pp_mw^*_q=
-(L^m|_{x/q})^*_{+}w^*_q
\end{equation}
One uses the standard notation ${\bf (P37)}\,\,res_z(\sum a_iz^i)=a_{-1}$ and 
$res_{D_q}(\sum b_iD^i_q)=b_{-1}$ and proves an analogue of Dickey's lemma, namely
\bq\label{13h}
res_z(Pexp_q(xz)Q^*|_{x/q}exp_{1/q}(-xz))=res_{D_q}(PQ)
\end{equation}
Further a q-bilinear identity is proved in the form ${\bf (P38)}\,\,res_z
(D^n_q\pp^{\ga}w_qw_q^*)=0$; the converse is also true.  In addition 
given formal series
\bq\label{14h}
w_q=\left(1+\sum_1^{\infty}w_iz^{-i}\right)exp_q(xz)exp\left(
\sum_1^{\infty}t_iz^i\right);
\end{equation}
$$w_q^*=\left(1+\sum_1^{\infty}w_i^*z^{-i}\right)exp_{1/q}(-xz)exp\left(
-\sum_1^{\infty}t_iz^i\right)$$
with {\bf (P38)} holding for any $n\in{\bf Z}_{+}$ and any multi-index $\ga$ with
nonnegative components $\ga_i$, then the operator $L=SD_qS^{-1}$ where $S=1+\sum
w_iD_q^{-i}$ is a solution of the q-KP hierarchy with wave and adjoint wave
functions given by $w_q$ and $w_q^*$.  As a consequence one can prove the existence
of a quantum tau function.  Indeed let ${\bf (P39)}\,\,\tl{w}_q=[exp_q(xz)]^{-1}
w_q$ and $\tl{w}_q^*=[exp_{1/q}(-xz)]^{-1}w^*_q$.  A little argument shows that
there is a function $\tau_q(x;t)$ such that
\bq\label{15h}
\tl{w}_q=\frac{\tau_q(x;t-[z^{-1}])}{\tau_q(x;t)}exp\left(\sum_1^{\infty}
t_iz^i\right);\,\,
\tl{w}_q^*=\frac{\tau_q(x;t+[z^{-1}])}{\tau_q(x;t)}exp\left(-\sum_1^{\infty}t_i
z^i\right)
\end{equation}
Equivalently then one can write
\bq\label{16h}
w_q=\frac{\tau_q(x;t-[z^{-1}])}{\tau_q(x;t)}exp_q(xz)exp\left
(\sum_1^{\infty}t_iz^i\right);
\end{equation}
$$w_q^*=\frac{\tau_q(x;t+[z^{-1}])}{\tau_q(x;t)}exp_{1/q}(-xz)exp\left(-\sum_1^{\infty}
t_iz^i\right)$$
It follows that if $L_1=\pp+\sum_1^{\infty}a_i\pp^{-i}$ with 
$\pp/\pp t_1$ is a solution of the KP hierarchy with tau function $\tau$, then
\bq\label{17h}
\tau_q(x;t)=\tau(t+[x]_q);\,\,[x]_q=\left(x,\frac{(1-q)^2}{2(1-q^2)}x^2,
\frac{(1-q)^3}{3(1-q^3)}x^3,\cdots,\right)
\end{equation}
is a tau function for the q-KP hierarchy.  Finally applications to N-qKdV are given
and in particular for $L=D_q^2+(q-1)xuD_q+u$ a solution for qKdV one has
\bq\label{18h}
u(x;t)=D_q\frac{\pp}{\pp t_1}log\tau_q(x;t)\tau_q(xq;t)
\end{equation}
This via $S=1-(\pp/\pp t_1)log\tau_q(x;t)D_q^{-1}+\cdots$ and
$S^{-1}=1+(\pp/\pp t_1)log\tau_q(x;t)D_q^{-1}+\cdots$.
\\[3mm]\indent
We recall also the standard symbol calculus for PSDO following e.g.
\cite{gam,mpa,tp}.
First one recalls from \cite{gam} the ring 
$\mf{A}$ of pseudodifferential
operators (PSDO) via PSD symbols (cf. also \cite{tp} for a more mathematical
discussion).  Thus one looks at formal series ${\bf
(P40)}\,\,A(x,\xi)=\sum_{-\infty}^na_i(x)\xi^i$ where $\xi$ is the symbol for
$\pp_x$ and $a_i(x)\in{\bf C}^{\infty}$ (say on the line or circle).  The
multiplication law is given via the Leibnitz rule for symbols ${\bf
(P41)}\,\,A(x,\xi)\ci B(x,\xi)=\sum_{k\geq 0}(1/k!)A^k_{\xi}(x,\xi)B_x^{(k)}(x,\xi)$
where $A_{\xi}^k(x,\xi)=\sum_{-\infty}^na_i(x)(\xi^i)^{(k)}$ and $B_x^{(k)}(x,\xi)=\sum_
{-\infty}^nb_i^{(k)}(x)\xi^i$ with $b_i^{(k)}(x)=\pp_x^kb_i(x)$.  This determines a
Lie algebra structure on $\mf{A}$ via ${\bf (P42)}\,\,[A,B]=A\ci B-B\ci A$.  Now let
A be a first order formal PSDO of the form ${\bf
(P43)}\,\,A=\pp_x+\sum_{-\infty}^{-1}a_i(\tl{x})\pp_x^i$ where $\tl{x}\sim (x,t_2,
t_3,\cdots)$.  Then the KP hierarchy can be written in the form ${\bf (P44)}\,\,(\pp
A/\pp t_m)=[(A^m)_{+},A]$ which is equivalent to a system of evolution equations
${\bf (P45)}\,\,(\pp a_i/\pp t_m)=f_i$ where the $f_i$ are certain universal
differential polyomials in the $a_i$, homogeneous of weight $m+|i|+1$ where
$a^j_{-i}$ has weight $|i|+j+1$ for $a^j\sim\pp_x^ja$.
\\[3mm]\indent
Somewhat more traditionally (following \cite{tp} - modulo notation and various
necessary analytical details), one can write
\bq\label{1h}
Au(x)=(2\pi)^{-1}\int e^{ix\cdot\xi}a(x,\xi)\hat{u}(\xi)d\xi
\end{equation}
where $\hat{u}(\xi)=\int exp(-ix\cdot \xi)u(x)dx$.  One takes $D=(1/i)\pp_x$ and
writes $a=symb(A)$ with $A=op(a)\sim\dot{A}$ where the $\cdot$ is to mod out
$\mf{S}^{-\infty}$ (we will not be fussy about this and will simply use A).  The
symbol of $A\ci B$ is then formally
\bq\label{2h}
(a\odot b)(x,\xi)=\sum\frac{1}{\ga!}\pp_{\xi}^{\ga}a(x,\xi)D_x^{\ga}b(x,\xi)
\end{equation}
corresponding to {\bf (P41)}, while $[A,B]=AB-BA$ corresponds to the symbol
${\bf (P46)}\,\,\{a,b\}=(\pp a/\pp\xi)(\pp b/\pp x)-(\pp a/\pp x)(\pp b/\pp \xi)$.
One notes that $\widehat{P(D)T}=P(\xi)\hat{T}$.
In any event it is clear that the algebra of differential operators on a manifold M
(quantum operators) may be considered as a noncommutative deformation of the algebra of
functions on $T^*M$ defined by canonical quantization via the symplectic form $\go=\sum
dp_i\wg dx^i$.  The extension to PSDO brings one into the aena of integrable systems
etc. Thus in a certain sense KP is an extension or generalization of quantum mechanics
(QM) based on the ring of PSDO (PSDO of all orders arise via $L^n_{+}$ in the higher
flows).
\\[3mm]\indent
In \cite{mpa} for example one extends matters to q-derivatives
$\pp_qf(z)=[f(qz)-f(z)]/ (q-1)z$ via ${\bf
(P47)}\,\,\pp_q(fg)=\pp_q(f)g+\tau(f)\pp_qg$ where
$\tau(f)(z)=f(qz)$ (note $\pp_q\tau=q\tau\pp_q$).  PSDO are defined via 
${\bf (P48)}\,\,A(x,\pp_q)=\sum_{-\infty}^nu_i(x)\pp^i_q$ with $\pp_qu=
(\pp_q u)+\tau(u)\pp_q$ and one has
\bq\label{6h}
\pp_q^{-1}u=\sum_{k\geq 0}(-1)^kq^{-k(k+1)/2}(\tau^{-k-1}(\pp_q^ku))\pp_q^{-k-1};
\end{equation}
$$\pp_q^nu=\sum_{k\geq 0}\left[\begin{array}{c}
n\\
k\end{array}\right]_q(\tau^{n-k}(\pp_q^ku))\pp_q^{n-k}$$
Recall here 
\bq\label{7h}
(n)_q=\frac{q^n-1}{q-1};\,\,\left[\begin{array}{c}
m\\
k\end{array}\right]_q=\frac{(m)_q(m-1)_q\cdots(m-k+1)_q}{(1)_q(2)_q\cdots (k)_q}
\end{equation}
Then the q-analogue of the Leibnitz rule is
\bq\label{8h}
A(x,\pp_q)B(x,\pp_q)=\sum_{k\geq
0}\frac{1}{(k)_q!}\left(\frac{d^k}{d\pp_q^k}A\right)*(\pp_q^kB)
\end{equation}
$$\frac{d^k}{d\pp_q^k}(f\pp_q^{\ga})=(\ga)_q(\ga-1)_q\cdots
(\ga-k+1)_qf\pp_q^{\ga-k}$$
One also uses the rules ${\bf (P49)}\,\,f*\pp_q=f\pp_q,\,\,\pp_q*f=\tau(f)\pp_q$,
and
$\pp_q^{-1}*f=\tau^{-1}(f)\pp_q^{-1}$.  Then set
$L_q=\pp_q+u_1(z)+u_2(z)\pp_q^{-1}+ u_3(z)\pp_q^{-2}+\cdots$ and one has q-KP via
$(\pp L_q/\pp t_m)=[L_q,(L_q^m)_{+}]$ where the order is different in $[\,,\,]$ and
$u_1(z)$ has a nontrivial evolution because of {\bf (P49)}.
\\[3mm]\indent
In accord with the procedures of \cite{gal,sh} we should now represent the ring
$\mf{A}_q$ of qPSDO symbols via a product as in say {\bf (P41)} and thence provide 
expressions for deformation thereof.  The $X$ and $P$ variables should come from the
phase space for dKP.  Evidently the qPSDO symbols will involve a variation on
\eqref{8h} (cf.
\eqref{2h}) and in view of the lovely development sketched above from \cite{ia} it
should be better to phrase matters in that notation.  Thus use {\bf (P33)} where 
$D_q\sim\pp_q$ and from {\bf (P28)} one has $D_q(fg)=f(qx)D_qg+(D_qf)g(x)$.  Thus the
rules of \cite{mpa} should apply to $D_q$ with suitable embellishment and we look at
(cf. {\bf (P29)} and {\bf (P30)})
\bq\label{19h}
\sum_{-\infty}^na_i(x)D_q^i\sum_{-\infty}^nb_j(x)D^j_q=\sum_{-\infty}^na_i(x)\sum_{-\infty}^n
\sum_0^{\infty}\left(\begin{array}{c}
i\\
k\end{array}\right)_q(D^{i-k}D_q^kb_j(x))D_q^{i+j-k}=
\end{equation}
$$=\sum_{k\geq 0}\sum_{i,j} a_i(x)\left(\begin{array}{c}
i\\
k\end{array}\right)_q(D^{i-k}D_q^kb_j(x))D_q^{i+j-k}$$
This amounts to ($\xi\sim D_q$)
\bq\label{20h}
\sum a_i(x)\xi^i\ci\sum b_j(x)\xi^j=\sum_{i,j,k}a_i(x)\left(\begin{array}{c}
i\\
k\end{array}\right)_q(D^{i-k}D_q^kb_j(x))\xi^{i+j-k}
\end{equation}
leading to
\bq\label{22h}
\{a,b\}=\sum_{k\geq
0}\sum_{i,j}\xi^{i+j-k}\left[a_i(x)
\binom{i}{k}_qD^{i-k}D^k_qb_j(x)-b_j(x)\binom{j}{k}_qD^{j-k}D_q^ka_i(x)
\right]
\end{equation}
Another way of writing this could be based on (cf. {\bf (P41)})
\bq\label{23h}
\frac{1}{(k)_q!}A_{\xi}^k\sim\sum_i\left(\begin{array}{c}
i\\
k\end{array}\right)_qa_i(x)\xi^{i-k}D^{i-k};\,\,B_x^{(k)}=\sum_jD^k_qb_j(x)\xi^j
\end{equation}
thus as symbols
\bq\label{24h}
a\ci b\sim\sum_{k\geq 0}\frac{1}{(k)_q!}A_{\xi}^kB^{(k)}_x;\,\,\{a,b\}\sim
\sum_{k\geq 0}\frac{1}{(k)_q!}\left(A^k_{\xi}B^{(k)}_x-B_{\xi}^kA_x^{(k)}\right)
\end{equation}
more in keeping with \eqref{2h}.  In this direction one could write e.g.
$\eta=\gk^{-1}\xi$ with
$\pp_{\eta}=\gk\pp_{\xi}$ and ${\bf
(P50)}\,\,(1/(k)_q!)\pp_{\eta}^kA\sim(\gk^k/(k)_q!)\pp_{\xi}^kA$ (cf. \eqref{23h}). 
Then define ${\bf (P51)}\,\,\{a,b\}_{\gk}=\sum_{k\geq
0}(\gk^k/(k)_q!)(A_{\xi}^kB^{(k)}_x- B^k_{\xi}A_x^{(k)})$.
In any event we have shown heuristically (see Section 4 for more detail
and enhancement)
\begin{proposition}
The calculi of PSDO and q-PSDO correspond symbolically via $\pp\sim D_q=\pp_q$
and suitable insertion of $D\sim \tau$ factors along with q-subscripts
(individual terms may differ because e.g. brackets
$[\,\,,\,\,]$ have different degrees, etc.).
\end{proposition}
\indent
We note from \cite{dam} an associative q-Weyl type star product based on $(\nu=log(q)$
and $\hbar=0$)
\bq\label{25h}
*^q_W=exp\left(-\frac{\nu}{2}\left[\ola{\pp}_xxp\ora{\pp}_p-\ola{\pp}_ppx\ora{\pp}_x
\right]\right)
\end{equation}
which essentially corresponds to our formula (4.14) in \cite{czj}, namely
\bq\label{26h}
\{f,g\}_M\sim f(xq^{-p\pp_p},pq^{x\pp_x})g(x,p)-g(xq^{-p\pp_p},pq^{x\pp_x})f(x,p)
\end{equation}
This is however quite different from our version \eqref{2t} 
(with $X^D_fg\sim\{f,g\}_M$) and it
is different also from \eqref{58}.  Some obstacles to the use of 
a $D_q$ version of \eqref{26h} for a Weyl ordered q-plane with $PX-qXP=i\hbar$ are
discussed in \cite{dam}.  In particular the lack of a complete basis is indicated
and it may be that some version of
\eqref{58} will circumvent this problem.  We note also
that {\bf (P51)} for example arising from the PSDO calculus \eqref{19h} - \eqref{24h}
etc. is different from the associative star products suggested in \cite{dam} of the form
$*_W^q$ above and ($\nu=log(q)$ and $\hbar=0$)
\bq\label{27h}
*^q_S=exp(\nu\ola{\pp}_ppx\ora{\pp}_x);\,\,*^q_A=exp(-\nu\ola{\pp}_xxp\ora{\pp}_p)
\end{equation}
Here we apparently must identify \eqref{24h} and {\bf (P51)} as the algebra of qPSD
symbols or operators while expressions such as \eqref{25h} or \eqref{27h} (and
perhaps \eqref{2t}) would correspond to q-Moyal type deformations of
dKP; their relations should then be determined.
We note also the (possibly nonassociative)
$D_q$ versions of $*^q_S$ and $*^q_A$ from \cite{dam}
for $\hbar\ne 0$, namely (here apparently $[m]\sim (m)_q\sim (q^m-1)/(q-1)$)
\bq\label{28h}
*_S=\sum_{r=0}^{\infty}\frac{(i\hbar)^r}{(r)_q!}\ola{D}_p^rexp(\nu\ola{\pp}_ppx\ora
{\pp}_x)\ora{D}^r_x;
\end{equation}
$$*_A=\sum_{s=0}^{\infty}(-\nu\ola{\pp}_xx)^s\sum_{r=0}^{\infty}\frac{(-i\hbar)^r
q^{r(r-1)/2}}{(r)_q!}\ola{D}_x^r\ora{D}^r_p(p\ora{\pp}_p)^s$$
We recall from \cite{dzz} that standard ordering S here means
$x^mp^n\to\hat{x}^m\hat{p}^n$ while antistandard A means $x^mp^n\to
\hat{p}^n\hat{x}^m$.  Thus we want to compare {\bf (P51)} with \eqref{25h}, or 
\eqref{27h} or \eqref{28h} and in view of the $D_q$ operators we probably want
\eqref{28h} or a Weyl form of this.  Recall that for
$\nu=log(q)$ one has $exp(\nu p\pp_p)f(p)=f(qp)$ and e.g.
\bq\label{29h}
f(pq^{-x\pp_x})g(x)\sim \sum f_np^nq^{-nx\pp_x}g(x)=\sum f_np^ng(q^{-n}x)
\end{equation}
Thus writing out the first equation in \eqref{28h} for example one gets
($\mf{D}=D_q$)
\bq\label{30h}
f*_Sg=\sum_0^{\infty}\frac{i\hbar)^r}{(r)_q!}\mf{D}_p^rf(x,p)e^{\nu\ola{\pp}_ppx\ora{\pp}_x}
\mf{D}_x^rg(x,p)
\end{equation}
\indent
Before embarking on questions of comparison \`a la \cite{gal,sh} let us recall some
results from \cite{fai}.  Here one considers star products of the form
${\bf (P52)}\,\,f\bgs g=fg+\sum_{n\geq 1}h^nB_n(f,g)$ where the
$B_n$ are bilinear differential operators. In particular in \cite{fai} one shows that
any bracket of the form
\bq\label{n9}
\{f,g\}=\sum_{r=1}^{\infty}\sum_{s=1}^{\infty}\gl^{r+s-2}\sum_{j=0}^r\sum_{k=0}^s
b_{rj,sk}(\pp_x^j\pp_y^{r-j}f)(\pp_x^k\pp_y^{s-k}g)
\end{equation}
may be transformed to one with $b_{00,10}=b_{00,11}=0$ and any such bracket satisfying
the Jacobi identity must be of the form
\bq\label{n10}
\{f,g\}=\sum_{r=1}^{\infty}\gl^{r-1}\sum_{j=0}^r\sum_{k=0}^sb_{rjk}(\pp_x^j\pp_y^{r-j}
f)(\pp_x^k\pp_y^{r-k}g)
\end{equation}
By suitable hocus pocus one shows also that \eqref{n10} plus Jacobi is equivalent to
Moyal.  Note here that \eqref{n18} below is of this form with $b_{rr0}\ne 0,\,\,
b_{r0r}\ne 0$, and all other coefficients equal 0.  Also $b_{110}=b_{101}$ (as required)
and the Jacobi condition for $\{f,g\}=(1/h)(f\bgs g-g\bgs f)$ can be proved directly via
associativity of $\bgs$ (exercise).  
Thus
\bq\label{n11}
\{\{f,g\},h\}+\{\{h,f\},g\}+\{\{g,h\},f\}=0
\end{equation}
\\[3mm]\indent
{\bf REMARK 3.1.}
We recall here the argument in \cite{sh} relating KP and $dKP_M$.  Here
\bq\label{n12}
\{f,g\}_{\gk}=\sum_0^{\infty}\frac{(-1)^s\gk^{2s}}{(2s+1)!}\sum_{j=0}^{2s+1}
(-1)^j\binom{2s+1}{j}(\pp_x^j\pp_p^{2s-j+1}f)(\pp_x^{2s-j+1}\pp_p^jg)
\end{equation}
gives an expression for $dKP_M.\hfill\bs$
\\[3mm]\indent
In order to establish an equivalence to the KP hierarchy, based on PSDO of the form
\bq\label{n13}
\pp_nL=[(L^n)_{+},L];\,\,L=\pp+\sum_2^{\infty}v_n(x,t_i)\pp^{-n}
\end{equation}
one looks at
\bq\label{n14}
\pp_n\gl=\{(\mf{L}^n)_{+},\mf{L}\};\,\,\mf{L}=p+\sum_2^{\infty}u_n(x,t_i)p^{-n+1}
\end{equation}
and compares \eqref{n12} with a bracket allegedly based on PSDO of the form
(cf. \cite{sh})
\bq\label{n15}
\{f,g\}'_{\gk}=\sum_0^{\infty}\frac{\gk^{2n+1}}{(2n+1)!}\left[\pp_{\xi}^{2n+1}f\pp_x^
{2n+1}g-\pp_{\xi}^{2n+1}g\pp_x^{2n+1}f\right]
\end{equation}
Note here from {\bf (P41)}, \eqref{1h}, etc. 
${\bf (P53)}\,\,A\ci B=\sum(1/k!)\pp_{\xi}^kA(x,\xi)\pp_x^kB(x,\xi)$ so 
(cf. {\bf (P41)} and \cite{ch})
\bq\label{n16}
A\ci_{\gk}B=Ae^{\gk\ola{\pp}_{\xi}\ora{\pp}_x}B=\sum\frac{\gk^n}{n!}\pp_{\xi}^nA
\pp_x^nB
\end{equation} 
and the bracket based on this is not obviously the same as \eqref{n15} or equivalent.
Note also
\bq\label{n17}
A\ci_{\gk}B=A(x, \xi+\gk\pp_x)B(x,\xi);\,\,B\ci_{\gk}A=B(x,\xi+\gk\pp_x)A(x,\xi)
\end{equation}
We see however that ${\bf (P54)}\,\,(1/\gk)A\ci_{\gk}B-B\ci_{\gk}A)=\{A,B\}_{\gk}$ is of
the form \eqref{n10} with $b_{rr0}\ne 0,\,\,b_{r0r}\ne 0$, and all other coefficients
equal 0.  Also $b_{110}=-b_{101}$ and the Jacobi identity will follow from
associativity so in fact a bracket such as {\bf (P54)} is equivalent to Moyal in the
symbols involved.  Note here that associativity is not obvious however (although
asserted in \cite{kt,mpb}) and \eqref{n18} below was only asserted to be associative in
\cite {mzm} (not proved).  Let us clarify this since it is not entirely trivial.
Thus consider \eqref{n16} along with e.g. (more on this below - cf. Remark 3.3)
\bq\label{n18}
f\bgs g=fg+\sum_{n\geq 1}\frac{(-h)^n}{n!}p^n\pp_p^nfx^n\pp_x^ng
\sim [q^{-x'\pp_{x'}p\pp_p}f(x,p)g(x',p')]|_{x,p} 
\end{equation}
while ${\bf
(P55)}\,\, A\ci_{\gk}B\sim e^{\gk\pp_{xi}\pp_{x'}}A(x,\xi)B(\xi',x')\,\,(x'\to
x,\,\xi'\to \xi)$. In fact associativity for $A\ci_{\gk}B$ is proved in \cite{kku}
(note the $\gk$ can be absorbed in $\xi$ by rescaling).  The trick is to use the formula
\bq\label{n20}
(a\xi^n)(b\xi^r)=a\sum_{k\geq 0}\binom{n}{k}\pp_x^kb\xi^{n-k+r}
\end{equation}
which shows that (for $\gk=1$)
\bq\label{n21}
A\ci_{\gk}B=\sum\frac{1}{m!}\pp_{\xi}^mA\pp_x^nB=
\end{equation}
$$=\sum\frac{1}{m!}\sum
a_nn(n-1)\cdots(n-m+1)\xi^{n-m}\cdot\sum b_j^{(m)}\xi^j=\sum a_n\binom{n}{m}b_j^
{(m)}\xi^{n-m+j}$$
Now for associativity one checks that ${\bf
(P56)}\,\,[\xi^n(a\xi^r)]b=\xi^n[(a\xi^n)b]$.  Thus the left side of {\bf (P56)} is 
\bq\label{n22}
\left[\sum_{\gag\geq 0}\binom{n}{\gag}a^{(\gag)}\xi^{n+r-\gag}\right]b=
\sum_{\gag,\mu\geq 0}a^{(\gag)}b^{(\mu)}\binom{n}{\gag}\binom{n+r-\gag}{\mu}\xi^
{n+r-\gag-\mu}
\end{equation}
and the right side of {\bf (P56)} is
\bq\label{n23}
\xi^n\left[a\sum_{\gag\geq 0}\binom{r}{\ga}b^{(\ga)}\xi^{r-\ga}\right]=
\sum_{\ga,\gb\geq 0}(ab^{(\ga)})^{(\gb)}\binom{n}{\gb}\binom{r}{\ga}\xi^{n+r-\ga-\gb}=
\end{equation} 
$$=\sum_{\ga,\gb,\gag\geq 0}a^{(\ga)}b^{(\gb+\ga-\gag)}\binom{\gb}{\gag}\binom
{n}{\gb}\binom{r}{\ga}\xi^{n+r-\ga-\gb}$$
For {\bf (P56)} one needs then ($n,r\in{\bf Z},\,\,\gag,\mu\in {\bf Z}_{+}$)
\bq\label{n24}
\binom{n}{\gag}\binom{n+r-\gag}{\mu}=\sum_{\ga+\gb=\gag+\mu,\,\ga,\gb>0}
\binom{\gb}{\gag}\binom{n}{\gb}\binom{r}{\ga}
\end{equation}
To prove \eqref{n24} one can start with 
\bq\label{n25}
\gag!\binom{n}{\gag}(1+x)^{n-\gag}(1+y)^r|_{y=x}=\left.\left[\frac{1}{\gag!}\pp_x^{\gag}
(1+x)^n(1+y)^r\right]\right|_{y=x}
\end{equation}
and picking out the coefficients of $x^{\mu}$ on both sides gives
\bq\label{n26}
\binom{n}{\gag}\binom{n+r-\gag}{\mu}=\frac{1}{\gag!}\left[x^{\mu}-coeff\left\{
\frac{1}{\gag!}\pp_x^{\gag}\sum\binom{n}{\gb}x^{\gb}\binom{r}{\ga}y^{\ga}\right\}_
{y=x}\right]=
\end{equation} 
$$=x^{\mu}-coeff\left\{\binom{\gb}{\gag}\sum\binom{n}{\gb}\binom{r}{\ga}
x^{\gb-\gag}x^{\ga}\right\}=\sum_{\ga+\gb=\gag+\mu}\binom{\gb}{\gag}
\binom{n}{\gb}\binom{r}{\ga}$$
(cf. also \cite{man}).  Associativity for $\bgs$ conceivably follows along similar
lines.
\\[3mm]\indent
{\bf REMARK 3.2.}
Thus, although the origin of \eqref{n15} is unclear, it is sufficient to work from
the original PSDO bracket \eqref{n16} and use \cite{fai} to assert that all of these 
brackets are equivalent to Moyal (equivalent means up to a change of variables). 
A specific correspondence as in \cite{gal} is not
needed then to assert indirectly that q-KP is equivalent to $dKP_M$ as in \eqref{n18}
once it is clear that q-KP corresponds to KP.  $\hfill\bs$
\\[3mm]\indent
{\bf REMARK 3.3.}  The form \eqref{n18} comes from a q-plane construction as in
\cite{mzm} (cf. also \cite{wq} and Remark 3.5).  Thus one writes ${\bf
(P57)}\,\,\hat{x}\hat{p}= q\hat{p}\hat{x}$ and ${\mf A}_x={\bf
C}[\hat{x}_1,\hat{x}_2]/{\mf R}$ where
${\mf R}\sim {\bf (P57)}$.  
We recall that if $x^i\sim\hat{x}^i$ are corresponding
commuting variables (corresponding e.g. to some ordering and isomorphism as in \cite
{czd,czj,cze,oh}); then using a Fourier transform
\bq\label{n1}
\tl{f}(k)=\frac{1}{2\pi}\int d^2xe^{-ik_jx^j}f(x)
\end{equation}
a unique operator
\bq\label{n2}
W(f)=\frac{1}{2\pi}\int d^2ke^{ik_j\hat{x}^j}\tl{f}(k)
\end{equation}
replaces $x^i$ by $\hat{x}^i$ in the most symmetric manner (Weyl quantization).
If the $\hat{x}^i$ have Hermitian properties then $W(f)$ will inherit them for real 
$f$.  Operators defined by \eqref{n2} can be multiplied and one wants to associate them
with classical functions.  If such a function exists we call it $f\bgs g$ defined via
${\bf (P58)}\,\,W(f)W(g)=W(f\bgs g)$ or more explicitly
\bq\label{n3}
W(f)W(g)=\frac{1}{(2\pi)^2}\int d^2kd^2\ell e^{ik_i\hat{x}^i}e^{i\ell_j\hat{x}^j}
\tl{f}(k)\tl{g}(\ell)
\end{equation}
If the product of exponents can be defined via the Baker-Campbell-Hausdorff (BCH)
formula then $f\bgs g$ will exist.  This is the case for a canonical structure
${\bf (P59)}\,\,exp(ik_i\hat{x}^i)exp(i\ell_j\hat{x}^j)=exp(i(k_j+\ell_j)\hat{x}^j
-(i/2)k_i\ell_j\gt^{ij})$.  In fact one can compute $(f\bgs g)(x)$ from \eqref{n3}
and {\bf (P59)} by replacing $\hat{x}$ with x, i.e.   
\bq\label{n4}
f\bgs g=\frac{1}{(2\pi)^2}\int d^2kd^2\ell e^{(ik_j+\ell_j)x^j-(1/2)k_i\gt^{ij}\ell_j}
\tl{f}(k)\tl{g}(\ell)=
\end{equation}
$$=e^{(1/2)\pp_{x^i}\gt^{ij}\pp_{y^j}}f(x)g(y)|_{y\to x}$$
which is the Moyal product.
For the q-plane the BCH formula cannot be used explicitly and the Weyl quantization
\eqref{n2} does not seem to be the most natural one (cf. \cite{dam}).
For now, in terms of algebraic structure only, any unique prescription of an operator
with a function of classical variables will suffice.  This could be e.g. normal order
where $\hat{x}$ operators are placed to the left of $\hat{y}$ operators (or
better here $\hat{x}^m\hat{p}^n$ corresponds to $\hat{x}\hat{p}$ order).  Thus define
${\bf (P60)}\,\,W(f(x,p))=:f(\hat{x},\hat{p}):$ and {\bf (P58)} becomes
${\bf (P61)}\,\,:f(\hat{x},\hat{p})::g(\hat{x},\hat{p}):=:(f\bgs g)(\hat{x},\hat{p}):$
which for monomials is
\bq\label{n5}
\hat{x}^m\hat{p}^n\hat{x}^a\hat{p}^b=q^{-na}\hat{x}^{m+n}\hat{p}^{a+b};\,\,
:\hat{x}^n\hat{p}^n::\hat{x}^a\hat{p}^b:=q^{-na}:\hat{x}^{m+n}\hat{p}^{a+b}:=
\end{equation}
$$=\left.W\left(q^{-x'\pp_{x'}p\pp_p}x^mp^nx^{'a}p^{'b}\right|_{x'\to x;\,p'\to p}
\right)$$
This generalizes for power series to
\bq\label{n6}
f\bgs g=q^{-x'\pp_{x'}p\pp_p}f(x,p)g(x',p')|_{x'\to x;\,p'\to p}
\end{equation}
which is \eqref {n18}.
One could equally well have used $\hat{p}\hat{x}$ ordering or Weyl ordering here.
For mononomials of fixed degree the $\hat{x}\hat{p}$, or $\hat{p}\hat{x}$, or
Weyl ordered products form a basis.  We recall for $f=\sum f_{mn}x^mp^n$
\bq\label{nn6}
:f(\hat{x},\hat{p}):_W=\sum \frac {f_{mn}}{2^m}\sum_0^m\binom{m}{\ell}\hat{x}^{m-\ell}
\hat{p}^n\hat{x}^{\ell}
\end{equation}
For the q-plane the form {\bf (P58)} provides a formula ${\bf (P62)}\,\,
W(x^ix^j)=:\hat{x}^i\hat{x}^j:$ (with say $\hat{x}\hat{p}$ ordering) and 
leads to 
\bq\label{n8}
f\bgs g=fg +\sum_{n\geq 1}\frac{1}{n!}(-h)^n(p\pp_p)^nf(x\pp_x)^b
\end{equation}
which is \eqref{n18}.  Note for $f,g$ of the form $f=\sum f_{mn}x^mp^n$ with $m\geq 0$
and $-\infty<n\leq N$ \eqref{n18} or \eqref{n8} will have the same form.  
For this, working with $\hat{x}\hat{p}$ ordering, we recall 
$\hat{x}^m\hat{p}^n\hat{x}^a\hat{p}^b=q^{-na}\hat{x}^{m+a}\hat{p}^{n+b}$ and if
$n=-\eta$ we get $q^{\eta a}$ as a multiplier.  This is consistent with moving
$p^{-\eta}$ past $x^a$ with $p^{-1}x=qxp^{-1}$ (from $xp=qpx$).  Hence the formulas
\eqref{n5}, \eqref{n6}, and \eqref{n8} remain valid.
$\hfill\bs$
\\[3mm]\indent
In addition to the constructions for q-KP in \eqref{9h} - \eqref{17h} we mention here
the Frenkel (F) and Khesin-Lyubashenko-Roger (KLR) versions of q-KP (cf.
\cite{fa,kap}).  Thus write ($t=(t_1,t_2,\cdots),\,\,t_1\sim x$)
\bq\label{29s} 
Q=D+a_0(t)D^0+a_{-1}(t)D^{-1}+\cdots;\,\,Q_q=D_1+b_0(t)D_q^0+b_{-1}(t)D_q^{-1}+\cdots
\end{equation}
where $Df(x)=f(qx)$ and $D_qf(x)=[f(qx)-f(x)]/(q-1)x$.  The F and KLR hierarchies are
defined via
\bq\label{30s}
\frac{\pp Q}{\pp t_n}=[(Q^n)_{+},Q]\,\,\,(F);\,\,\frac{\pp Q_q}{\pp
t_n}=[(Q_q^n)_{+},Q_q]\,\,\,(KLR)
\end{equation}
and there is an isomorphism ${\bf (P63)}\,\,\hat{{}}:\,{\mf D}_q\to{\mf D}$ mapping the
F or KLR systems into the discrete KP hierarchy (cf. \cite{at,az,dan} - we omit here a
discussion of discrete KP).  These systems are equivalent by virtue of a correspondence
\bq\label{38s}
a_i(y)=\sum_{0\leq k\leq n-i}\frac{\left[\begin{array}{c}
k+i\\
k\end{array}\right]}{(-y(q-1)q^i)^k}b_{k+i}(y)
\end{equation}
Consider now a suitable space of functions $f(x)$ represented by ``Fourier" series ${\bf
(P64)}\,\,f(x)=\sum_{-\infty}^{\infty}f_n\phi_n(x)$ for $\phi_n(x)=\gd(q^{-n}y^{-1}x)$
for $q\ne 1$ and $y\in {\bf R}$.  Set $\gl_i=D^i\gl_0=\gl(yq^i)$; then the Fourier
transform $f\to{\mf F}f=(\cdots,f_n,\cdots)_{n\in{\bf Z}}$ induces an algebra
isomorphism
$\hat{{}}:\,{\mf D}_q\to{\mf D}$ via
\bq\label{36s}
\sum a_i(y)D^i\to\sum\hat{a}_i\gL^i=\sum diag(\cdots,a_i(xq^n),\cdots)_{n\in{\bf
Z}}\gL^i
\end{equation}
where $\gL\sim(\gd_{i,j-1})_{i,j\in{\bf Z}}$ is a shift operator.  In addition
\bq\label{37s}
\sum_0^nb_i(y)D^i_q=\sum_0^na_i(y)(-\gl
D)^i\to\gep\left(\sum_0^n\hat{a}_i\gL^i\right)\gep^{-1}
\end{equation} 
where
\bq\label{322s}
\gep=diag\left(\cdots,\gl_{-2}\gl_{-1},-\gl_{-1},1,-\frac{1}{\gl_0},
\frac{1}{\gl_0\gl_1},-\frac{1}{\gl_0\gl_1\gl_2},\cdots\right);\,\,\gep_0=1
\end{equation}
\indent
In any event we have seen that the calculation of PSDO  symbolically corresponds
symbolically to  qPSDO (Proposition 3.1).  Further from \cite{fai} we have seen that all
Moyal brackets based on an associative star product of the form {\bf (P52)} are
equivalent.   Consequently one has heuristically (cf. also Section 4)
\begin{proposition}
Given the associativity of \eqref{n18} asserted in \cite{mzm} it yields a Moyal bracket
equivalent to the standard one; this applies also to the associative star products
\eqref{25h} and \eqref{27h} and to \eqref{58} provided it is associative.  Hence
$dKP_M$ based on any such star product is equivalent to KP via \cite{gal,sh} and thence
to qKP symbolically as in Proposition 3.1 (embellished as in Proposition 4.1). 
The fact that q-plane structure led to \eqref{n18} is immaterial
here given its associativity.
\end{proposition}
\indent
This seems reasonable since knowing that KP is a ``quantization" of dKP under a
suitable Moyal bracket (and hence a generalized quantum theory)
one would expect the isomorphic theory qKP to be some kind of quantization of $dKP_M$
under equivalent Moyal brackets as indicated (cf. Proposition 3.1 and Proposition
4.1 for the correspondence $KP\leftrightarrow qKP$).
Explicit examples of comparisons as in
\cite{gal,sh} would be useful (cf. Remark 3.4).
\\[3mm]\indent
{\bf REMARK 3.4.}  For a start in this direction explicit calculations to compare e.g.
qKP and
$dKP_{\bgs}$ can be carried out using the compatibility conditions ${\bf
(P65)}\,\,D_xp=q^{-1}pD_x$ and
$D_px=qxD_p$ on the q-plane ($D\sim D_q$).  This can be confirmed as follows.
 From \cite{ka} one has
\bq\label{n28}
D_q(f(x)g(x))=g(x)D_qf(x)+f(qx)D_qg(x)
\end{equation}
Next note, in analogy to the formula from \cite{ka} (based on \eqref{1g}, \eqref{2g},
etc.)
\bq\label{n29}
\pp_i^q(f_N(x_N)\cdots f_1(x_1))=f_N(qx_N)\cdots
f_{i+1}(qx_{i+1})D_{q^2}f_i(x_i)f_{i-1}(x_{i-1})\cdots f_1(x_1)
\end{equation}
we have e.g.
\bq\label{n30}
D_p(f(x)g(p))=\frac{p^{-1}}{q-1}f(x)[g(qp)-g(p)]
\end{equation}
Now $xp=qpx$ so $p^{-1}x=qxp^{-1}$ and $p^{-1}x^m=q^mx^mp^{-1}$ leading to
\bq\label{n31}
D_p(f(x)g(p))=f(qx)\left[\frac{g(qp)-g(p)}{p(q-1)}\right]=f(qx)D_pg(p)
\end{equation}
Similarly ${\bf (P66)}\,\,D_x(f(x)g(p))=(D_qf(x))g(p)$ and we ask what this means for
$xD_p$ and $pD_x$.  From \eqref{n31} one has ${\bf (P67)}\,\,D_pf(x)=f(qx)D_p$ so 
$D_px=qxD_p$.  As in \eqref{n31} we look at $D_x(g(p)f(x))$ and note that 
$p^{-1}x=qx^{-1}p$ or $x^{-1}p=q^{-1}px^{-1}$ which means that $x^{-1}p^n=
q^{-n}p^nx^{-1}$ so
\bq\label{n32}
D_x(g(p)f(x))=g(q^{-1}p)\left[\frac{f(qx)-f(x)}{x(q-1)}\right]=g(q^{-1}p)D_xf(x)
\end{equation}
Consequently ${\bf (P68)}\,\,D_xg(p)=g(q^{-1}p)D_x$ and $D_xp=q^{-1}pD_x$ in agreement
with {\bf (P65)}. For the calculations one starts with ${\bf (P69)}\,\,
\gl=\mf{L}=p+a_0+\sum_1^{\infty}a_ip^{-i};\,\,\gl^2_{+}=\mf{L}^2_{+}=p^2+u_1p+
u_0$ where $a_i=a_i(x,t)$ and $u_i=u_i(x,t)$ with e.g. $\gl^2\sim\gl\bgs\gl$.  One
computes e.g. $\pp_1\gl=\{\gl^2_{+},\gl\}_{\bgs}=\gl_{+}^2\bgs\gl-\gl\bgs\gl^2_{+}$
and compares with $\pp_2L=[L^2_{+},L]$ based on $L=D_q=a_0+\sum a_iD_q^{-i}$ and
$L^2_{+}= D_q^2+u_1D_q+u_0$.  Under certain conditions ($a_0,a_1,u_0,u_1$ constant in 
$x$) compatible equations seem to appear, modulo solution of second order difference
equations
\bq\label{n56}
(D_q^2+2c_1D_q+c_2)a_n=2c_1(c_1+a_2+\cdots +a_{n-1})\,\,\,(n\geq 2)
\end{equation}
(cf. \cite{cze}).$\hfill\bs$
\\[3mm]\indent
{\bf REMARK 3.5.}
We note here in passing another way of dealing with q-plane differential operators
following \cite{oh} (cf. also \cite{czd,czj,cze}).  Thus one imbues
the q-plane 
or Manin plane (cf. \cite{mzt}) with the
natural associated covariant calculus (cf. \eqref{1g} - \eqref{2g} below).  Thus we
can treat the phase space as a q-plane
$x^1=x$ and
$x^2=p$ for say
$GL_q(2)$ with rules ($q\sim exp(\hbar)$)
\bq\label{1g}
xp=qpx;\,\,\pp_ix^j=qx^j\pp_i\,\,(i\ne
j);\,\,\pp_ix^i=1+q^2x^i\pp_i+q\gl\sum_{j>i}x^j\pp_j
\end{equation}
where $\gl=q-q^{-1}$.  Thus
\bq\label{2g}
\pp_xx=1+q^2x\pp_x+(q^2-1)p\pp_p;\,\,\pp_pp=1+q^2p\pp_p;\,\,\pp_xp=qp\pp_x;\,\,
\pp_px=qx\pp_p
\end{equation}
Note here for example
\bq\label{2gg}
\pp_pxp=qx\pp_pp=qx(1+q^2p\pp_p)=qx+q^3xp\pp_p;
\end{equation}
$$\pp_pxp=\pp_pqpx=q(1+q^2p\pp_p)x
=qx+q^3xp\pp_p$$
but a situation $p\sim \pp_x$ with $px-qxp=i\hbar$ is excluded (cf. \cite{dam}).
One denotes by $Diff_{q^2}(1)$ the ring generated by $x,\pp^q$ obeying ${\bf (P70)}\,\,
\pp^qx=1+q^2x\pp^q$ with $\pp^qf(x)=[f(q^2x)-f(x)]/(q^2-1)x=D_{q^2}f(x)$.  
We distinguish now scrupulously between $\pp^q\sim D_{q^2},\,\,D_q,$ and $\pp_x,\pp_p$
as normal q-derivatives. Introduce
${\bf (P71)}\,\,\mu_k=1+q\gl\sum_{j\geq k}x^j\pp_j$ so the last equation in \eqref{1g}
takes the form ${\bf (P72)}\,\,\pp_ix^i=\mu_i+x^i\pp_i$ (note the $\mu_i$ are
operators). Now there results
\bq\label{3g}
\mu_ix^j=x^j\mu_i\,\,(i>j)\,\,\mu_ix^j=q^2x^j\mu_i\,\,(i\leq
j);
\end{equation}
$$\mu_i\pp_i=\pp_j\mu_i\,\,(i>j)\,\,\mu_i\pp_j=q^{-2}\pp_j\mu_i\,\,(i\leq j)$$
which implies ${\bf (P73)}\,\,\mu_i\mu_j=\mu_j\mu_i$.  Next define ${\bf (P74)}\,\,X^i
=(\mu_i)^{-1/2}x^i$ and $D_i=q(\mu_i)^{-1/2}\pp_i$ from which follows
\bq\label{4g}
X^iX^j=X^jX^i;\,\,D_iD_j=D_jD_i;\,\,D_iX^j=X^jD_i\,\,(i\ne j);\,\,D_jX^j=1+q^{-2}X^jD_j
\end{equation}
Thus the relations in \eqref{1g} are completely untangled.  The $D_j$ correspond to
$D_{q^{-2}}$ and evidently $Diff_{q^2}(1)$ is isomorphic to $Diff_{q^{-2}}(1)$ since,
for $\gd^q=q\mu^{-1/2}\pp^q$ and $y=\mu^{-1/2}x$ with operators $x,\pp^q$ satisfying
{\bf (P70)}, one has $\gd^qy=1+q^{-2}y\gd^q$.  Further the ring isomorphism between
$Diff_{q^2}(1)$ (generated by ($x,\pp^q$) and $Diff(1)$ (generated by ($x,\pp$) 
can be established via e.g. ${\bf (P75)}\,\,\pp^q=(exp(2\hbar x\pp)-1)/x(q^2-1)$
(cf. {\bf (P70)}).  Thus $exp(2\hbar x\pp)-1=x(q^2-1)\pp^q$ or $2\hbar
x\pp=log[1+x(q^2-1)\pp^q]$.  Since the ring properties are not immediate from this one
can go to an alternative $\ul{noncanonical}$ isomorphism as follows (cf. \cite{oh}).
Let $x_c^i$ be classical commuting variables (here $x_c^1\sim x$ and $x_c^2\sim p$). 
Now choose some ordering of the nonclassical $x^i$ (e.g. Weyl ordering, or $xp$
ordering, or
$px$ ordering).  Then any polynomial $P(x)$ can be written in ordered form and
replacing $x^i$ by $x_c^i$ one gets a polynomial symbol $\gs(P)$ of classical variables
$x_c^i$.  This determines a symbol map $\gs:\,{\bf C}[x^i]\to{\bf C}[x_c^i]$ which is 
a noncanonical isomorphism (dependent on the choice of ordering)
between polynomial rings.  Then for any polynomial $\phi
(x_c^i)$ and any q-differential operator $D$ one writes ${\bf
(P76)}\,\,\hat{D}\phi=\gs(D(\gs^{-1}(\phi)))$, i.e. $\hat{D}$ is the composition
\bq\label{5g}
{\bf C}[x_c^i]\stackrel{\gs^{-1}}{\to}{\bf C}[x^i]\stackrel{D}{\to}{\bf C}[x^i]
\stackrel{\gs}{\to}{\bf C}[x_c^i]
\end{equation}
This provides a ring isomorphism of q-differential operators and classical differential
operators, the latter corresponding to polynomials in $(x,p,\pp_x,\pp_p)$ with 
relations $\pp_xx=x\pp_x+1,\,\,\pp_pp=p\pp_p+1,\,\,x\pp_p=\pp_px,$ and
$p\pp_x=\pp_xp$.  The explicit formulas will depend on the ordering and are determined
by $\hat{\pp}_i$ and $\hat{x}^i$.  Note ${\bf
(P77)}\,\,\widehat{D_1D_2}=\hat{D}_1\hat{D}_2$ since the $\hat{D}_i:\,\,{\bf
C}[x_c^i]\to{\bf C}[x_c^i]$ compose multiplicatively along with the $D_i:\,\,{\bf
C}[x^i]\to{\bf C}[x^i]$ under the given ordering.  To see this note from $\hat{D}_2\phi=
\gs(D_2(\gs^{-1}(\phi)))$ results
\bq\label{6g}
\hat{D}_1(\hat{D}_2\phi)=\gs(D_1(\gs^{-1}(\gs(D_2(\gs^{-1}(\phi))))))=\gs(D_1D_2(\gs^{-1}(\phi)))
\end{equation}
As for orderings, matters are clear for $xp$ or $px$ ordering and hence will also hold
for the completely symmetric Weyl ordering
\bq\label{7g}
x^np^m\sim\frac{1}{2^n}\sum_0^n\binom{n}{\ell}\tl{x}^{n-\ell}\tl{p}^m\tl{x}^{\ell}
\end{equation}
(cf. \cite{ch}).$\hfill\bs$

\section{DISCRETIZATION AND QUANTUM MECHANICS.}
\renewcommand{\theequation}{4.\arabic{equation}}
\setcounter{equation}{0}

We go back to earlier comments in Section 3 (after \eqref{2h}) which present KP as some
sort of extension of QM related via Moyal to a phase space $(X,P)$ and corresponding 
to dispersionless operators for dKP.  The formation of q-KP and various q-Moyal brackets
for dKP are present via Proposition 3.2 (in a somewhat dismissive manner)
and we want to examine this further.  Thus first look at QM and a q-QM obtained by
replacing differential operators in $\pp_i$ by q-difference operators using $\pp_i^q$
for example.  Take a 1-dimensional situation with $x,p\sim\pp_x$ as basic and recall
from \eqref{6h} that $\pp_qu=(\pp_qu)+\tau u\pp_q$ so the commutator relation
$[\pp_x,x]=1$ goes int ${\bf (P78)}\,\,[\pp_q,x]=\pp_qx-x\pp_q=\pp_qx+(\tau
x)\pp_q-x\pp_q=1+qx\pp_q-x\pp_q=1+(q-1)x\pp_q=1+(\tau-1)=\tau$.  More generally one has
Leibnitz formulas (cf. \eqref{6h})
\bq\label{y1}
\pp^nu=\sum_{k\geq 0}\binom{n}{k}\pp^ku\pp^{n-k}\longrightarrow \pp_q^nu=\sum_{k\geq
0}\binom{n}{k}_q\tau^{n-k}(\pp_q^ku)\pp_q^{n-k}
\end{equation}
so e.g. ${\bf (P79)}\,\,\pp^2u=(\pp^2u)+2\pp u\pp+u\pp^2\rightarrow
\tau^2u\pp_q^2+(1+q)\tau\pp_qu\pp_q+u\pp_q^2$.  Recall here 
\bq\label{yy1}
\binom{n}{k}_q=\frac{(q;q)_n}{(q;q)_k(q;q)_{n-k}};\,\,(q;q)_0=1;\,\,(q;q)_k=
\prod_1^k(1-q^j)
\end{equation}
Given that the algebra of differential operators ${\mc A}$ represents QM one can define
any isomorphic object or ``corresponding" object to also be a quantum theory (QT).  Thus
the correspondence $\pp\leftrightarrow\pp_q:\,\,{\mc A}\leftrightarrow {\mc A}_q$ with
Leibnitz rules as in \eqref{y1} leads to an algebra isomorphism.  Indeed we can simplify
here the formulas \eqref{19h} - \eqref{24h} via the symbol notation
\bq\label{y2}
\sum_0^na_i(x)D_q^i\sum b_j(x)D_q^j=\sum_{i=0}^na_i(x)\sum_{j=0}^n\sum_{k=0}^i
\binom{i}{k}_q\tau^{i-k}D_q^kb_j(x)D_q^{i-k+j}
\end{equation}
whereas
\bq\label{y3}
\sum_0^na_i\pp^i\sum_0^nb_j\pp^j=\sum_{i=0}^na_i(x)\sum_{k=0}^i\binom{i}{k}\pp^kb_j(x)
\pp^{i-k+j}
\end{equation}
Thus in symbol form ($\pp\sim\xi$)
\bq\label{y4}
\sum a_i\xi^i\ci\sum b_j\xi^j=\sum_{i,j,k}\frac{i(i-1)\cdots(i-k+1)}{k!}a_i\xi^{i-k}
\pp_x^kb_j(x)\xi^j=
\end{equation}
$$=\sum\binom{i}{k}\pp_{\xi}\xi^ia_i\pp_x^kb_j\xi^j=\sum_{k=0}\frac{1}{k!}
\pp_{\xi}^ka(x,\xi)\pp_x^kb(x,\xi)$$
and for $D_q\sim \xi$ (as in \eqref{20h} - \eqref{24h})
one takes $\pp_{\xi}\sim (\pp/\pp \xi)_q$ to produce
\bq\label{y5}
\sum a_i\xi^i\ci\sum b_j\xi^j=\sum_{k=0}^n\frac{1}{k_q!}\pp_{\xi}^ka\pp_q^kb;
\end{equation}
$$\frac{1}{k_q!}\pp_{\xi}^ka=\sum_ia_i\sum_{k=0}^i\binom{i}{k}_q\xi^{i-k}
\tau^{i-k}=\sum_ia_i\sum_{k=0}^i\frac{1}{k_q!}\left(\frac{\pp}
{\pp\xi}\right)^k_q\xi^i\tau^{i-k}$$
since $[i_q\cdots (i-k+1)_q/k_q!]=\binom{i}{k}_q$.  Now any product of two symbols has
the form $a\ci b$ as in \eqref{y4} and one wants to check uniqueness.  Thus suppose
\bq\label{y7}
L=\sum_{k=0}^n\frac{1}{k!}\pp_{\xi}^ka\pp_x^kb=\sum_{m=1}^n\frac{1}{m!}\pp_{\xi}^m
\hat{a}\pp_x^m\hat{b}=R
\end{equation}
We note that $\pp L/\pp b_0=a$ and $\pp R/\pp \hat{b}_0=\hat{a}$ so if we stipulate
that $\hat{b}_0=b_0$ then $a=\hat{a}$ (note that $b_0$ only appears undifferentiated
in the term $\pp_{\xi}^0a\pp_x^0b=a(\sum_0^nb_i\xi^i)$
and in any event $\pp(\pp^nb_0)/\pp b_0=\pp^n(1)=0$).  
Once we have $a=\hat{a}$ then
one can equate coefficients of $\pp_{\xi}^ka$ and $\pp_{\xi}^k\hat{a}$ (taking $m=k$)
to get $\pp_x^kb=\pp_x^k\hat{b}$ which will yield $b=\hat{b}$ (for analytic
$b,\hat{b}$).  Similar comments apply to \eqref{y5} and
consequently for $b,\hat{b}$ with equal $b_0$
we can look at ${\mc A}\leftrightarrow{\mc A}_q$ as an algebra map with
$a\ci b\leftrightarrow a_q\ci b_q\sim (a\ci b)_q$ where powers of $\tau$
must be inserted correctly (e.g. one could define here ${\bf (P80)}\,\,
(\pp_{\xi})^k_q\xi^i=
i_q\cdots (i-k+1)_q(\tau\xi)^{i-k}$).  On the other hand, if $b_0\ne \hat{b}_0$ set
$\tl{b}_0=\hat{b}_0+(b_0-\hat{b}_0)$ so that $L=\hat{R}\rightarrow L=\tl{R}-
\sum(1/m!)\pp_{\xi}^m\hat{a}\pp_x^m(b_0-\hat{b}_0)$ and $\tl{R}$ involves the same
$\tl{b}_0=b_0$.  Then $\pp L/\pp b_0=a=\pp\tl{R}/\pp b_0-\hat{a}=\hat{a} -\hat{a} =0$.
this says $L=R$ is possible for nontrivial $a$ only when $b_0=\hat{b}_0$.
Thus $a_q\ci b_q$ can be written as in \eqref{y5} with
$\pp_{\xi}$ taken as a q-derivative {\bf (P80)} in ${\mc A}_q$.
\\[3mm]\indent
Next we want to extend such arguments to the rings ${\mf A},\,{\mf A}_q$ of PSDO
and qPSDO as in \eqref{1h} - \eqref{24h}.  We can try the same procedure with
\bq\label{y8}
L=\sum_{k\geq 0}\frac{1}{k!}A_{\xi}^k(x,\xi)\pp^k_xB(x,\xi)=\sum_{k\geq 0}\frac{1}{k!}
\hat{A}_{\xi}^k(x,\xi)\pp_x^k\hat{B}(x,\xi)=R
\end{equation}
where $A_{\xi}^k=\sum_{-\infty}^na_i(x)\pp_{\xi}^k\xi^i$ and $\pp_x^kB(x,\xi)=
\sum_{-\infty}^n\pp_x^kb_j(x)\xi^j$.  The same reasoning applies for $b=\hat{b}_0$
producing $A=\hat{A}$ and $b_0\ne\hat{b}_0$ implies $A=0$.  For the q-derivatives one
modifies again the form of $(\pp_{\xi})_q^k$ as in {\bf (P80)} to obtain a ring
isomorphism ${\mf A}\leftrightarrow {\mf A}_q$.  Hence Proposition 3.1 can be improved
in the form
\begin{proposition}
The correspondence ${\mf A}\leftrightarrow {\mf A}_q$ of Proposition 3.1 can be viewed
as a ring isomorphism and thus one can claim that e.g. q-KP is a generalized QT.
\end{proposition}
\indent
{\bf REMARK 4.1.}
One must be careful in relating discretization and quantization (cf. \cite{ch} for
extensive comments on this).  For example the Moyal bracket can be obtained by taking
a continuous limit of a discrete dynamical bracket but discretization of the Moyal
bracket does not lead back to the discrete situation (cf. \cite{ch,ga}).
$\hfill\bs$
\\[3mm]\indent
{\bf REMARK 4.2.}
One can think of KP as an extended or generalized QT in two ways.  The first way
involving $x,\,\pp_x$ has already been indicated but one can equally well look at
$z,\,\pp_z$ or $L,M$ (cf. \cite{daz,day}).  This seems to be related to the development
of \cite{kaw,kas,kbe,kzq} indicated in Section 2, and a logical background is the idea
of z-operators and the action on the Grassmannian as in \cite{daz,day} (cf. also
\cite{mal,sz}).  We give here a little background.  Thus one writes z-operators in the
form ${\bf (P81)}\,\,G=G(\pp_z,z)=\sum_{j\geq 0,i\leq i_0}a_{ij}z^i\pp_z^j$ acting on
formal power series $f(z)=\sum_{-\infty}^{\infty}f_kz^k$ 
on the unit circle $S^1$ (convergence is not considered
here).  Operators ${\bf (P82)}\,\,G=1+\sum_{j\geq 0,i<0}a_{ij}z^i\pp_z^j$ are called
monic.  To any z-operator one assigns a PSDO ${\bf (P83)}\,\,{}^TG(x,\pp)=\sum a_{ij}
x^j\pp^i$ where T means $z\to\pp$ with $\pp_z\to x$ and the factors in
reverse order.  Evidently
\bq\label{y9}
f\cdot g\to{}^Tg\cdot{}^Tf;\,\,[z,\pp_z]=[x,\pp]=-1;
\end{equation}
$$G(\pp_z,z)e^{xz}={}^TG(x,\pp)
e^{xz}=G(x,z)e^{xz}={}^TG(x,z)e^{xz}$$
In particular this is an anti-isomorphism PSDO $\leftrightarrow$ z-operators.  One
recalls that the Grassmanian Gr consists of linear subspaces $V\subset H$
($H\sim\{f(z),\,z\in S^1\}$) such that the natural projection $V\to H_{+}$ is 1-1.
It is then well known (cf. \cite{cw,mal,sz}) that $V\in Gr\Rightarrow$ there exists 
a monic z-operator G such that $V=GH_{+}$ and if a z-operator preserves $H_{+}$ then it
involves only nonnegative powers of z.  Let now $t^*=(t_2,t_3,\cdots)$ and
$\xi^*=\sum_2^{\infty}t_kz^k$.  Then $V\in Gr\rightarrow Vexp(-\xi^*(t^*,z))\in Gr$ so
there is a monic z-operator $W(t^*,\pp_z,z)$ such that ${\bf (P84)}\,\,exp(-\xi^*)V=W
(t^*,\pp_z,z)H_{+}$.  Then as above ${\bf
(P85)}\,\,W(t^*,\pp_z,z)exp(xz)=W(t^*,x,\pp)exp(xz)=W(t^*,x,z)exp(xz)$.  Hence
$W(t^*,x,z)exp(xz)=
W(t^*,\pp_z,z)exp(xz)\in Vexp(-\xi(t^*,z))$ and $W(t^*,x,z)exp(\xi(t,z))\in V$
($\xi(t,z)=xz+\xi(t^*,z)$).  This means ${\bf (P86)}\,\,\psi_V(t,z)=W(t^*,x,z)
exp(\xi(t,z))$ is the Baker-Akhiezer (BA) function of the KP hierarchy based on
$L=W(t^*,x,\pp)\pp W^{-1}(t^*,x,\pp)$ related to $V\in Gr$.
\\[3mm]\indent
The flows $V\to G(t)V$ can be related to the Orlov-Schulman operators and the Virasoro
algebra in a natural manner (cf. \cite{daz,day}).  What one sees here is that the QT
features of KP based on $x,\,\pp$ pass directly via anti-isomorphism to QT features for
$\pp_z,\,z$ (or $M,\,L$).  In particular corresponding phase space variables could be
viewed as
$\gz,\,z$ (where $\gz\sim
\pp_z$ in the same way as $p\sim \pp$) and this gives us an entr\'ee to the use of 
$\gz,\,z$ as phase space variables as used in \cite{kaw} and in Section 2.  Thus
in an obvious way one can expect Moyal type theory based on phase space functions
$A(z,\gz)$ to be deformation quantization equivalent to the generalized QT of 
z-operators which in turn is anti-isomorphic to the KP theory.  The 
subsequent natural emergence of z-operators $z,\,\pp_z$ in various formulas then
allows one to formulate the dynamical theory directly in terms of Orlov-Schulman
operators as in Section 2 and this could all lead to further perspective on vertex
operators.
$\hfill\bs$


\begin{thebibliography}{ccc}




%
\bibitem{au} M. Adler, T. Shiota, and P. vanMoerbeke,
Phys. Lett. A, 194 (1994), 33-43; Comm. Math. Phys., 171 (1995), 547-588;
solv-int 9812006
%
\bibitem{al} M. Adler and P. vanMoerbeke,
Comm. Math. Phys., 147 (1992), 25-56
%
\bibitem{at} M. Adler and P. vanMoerbeke,
Comm. Math. Phys., 203 (1999),185-210 = solv-int 9912014
%
\bibitem{az} M. Adler, E. Horozov, and P. vanMoerbeke,
Phys. Lett. A, 242 (1998), 139-151
%
\bibitem{ch} R. Carroll,
Quantum theory, deformation, and integrability, North-Holland, 2000
%
\bibitem{czd} R. Carroll,
math.QA 0101072, Applicable Anal., to appear
%
\bibitem{czj} R. Carroll,
math.QA 0101229, Inter. Jour. Math. and Math. Sci., to appear
%
\bibitem{cw} R. Carroll,
Topics in soliton theory, North-Holland, 1991
%
\bibitem{cqq} R. Carroll,
Applic. Anal., 49 (1993), 1-31; 56 (1995), 147-164
%
\bibitem{cm} R. Carroll,
Nucl. Phys. B, 502 (1997), 561-593; Lect. Notes Physics 502, Springer, 1998, pp.
33-56
%
\bibitem{cn} R. Carroll,
Jour. Nonlin. Sci., 4 (1994), 519-544
%
\bibitem{co} R. Carroll and Y. Kodama,
Jour. Phys. A, 28 (1995), 6373-6378
%
\bibitem{cze} R. Carroll,
Q-analysis and integrability, part of a book in preparation
%
\bibitem{dam} O. Dayi,
Jour. Math. Phys., 39 (1998), 485-496
%
\bibitem{dan} L. Dickey,
Lett. Math. Phys., 48 (1999), 277-289
%
\bibitem{dap} L. Dickey,
Acta Appl. Math., 47 (1997), 243-321
%
\bibitem{daq} L. Dickey,
Soliton equations and Hamiltonian systems, World Scientific, 1991
%
\bibitem{daz} L. Dickey,
hep-th 9210155 and 9312015
%
\bibitem{day} L. Dickey,
Mod. Phys. Lett. A, 8 (1993), 1259-1272 and 1357-1377
%
\bibitem{dzz} G. Dunne,
Jour. Phys. A, 21 (1988), 2321-2335
%
\bibitem{fai} P. Fletcher,
Phys. Lett. B, 248 (1990), 323-328
%
\bibitem{fam} P. diFrancesco, P. Mathieu, and D. S\'en\'echal,
Conformal field theory, Springer, 1997; hep-th 9112063
%
\bibitem{fal} P. diFrancesco and P. Mathieu,
hep-th 9109042
%
\bibitem{fa} E. Frenkel,
IMRN, 2 (1996), 55-76
%
\bibitem{faj} M. Freeman and P. West,
hep-th 9208013
%
\bibitem{ga} D. Galetti and M. Ruzzi,
Physica A, 26 (1999), 473-491; 33 (2000), 1065-1082
%
\bibitem{gal} J. Gawrylczyk,
Jour. Phys. A, 28 (1995), 593-605
%
\bibitem{gam} J. Gawrylczyk,
Jour. Math. Phys., 36 (1995), 3461-3478
%
\bibitem{hy} L. Haine and P. Iliev,
Jour. Phys. A, 30 (1997), 7217-7227; Int. Math. Res. Notices, 6 (2000), 281-323
%
\bibitem{ia} P. Iliev,
Jour. Phys. A, 31 (1998), L241-L244; Lett. Math. Phys., 44 (1998), 187-200
%
\bibitem{kaw} R. Kemmoku,
Jour. Phys. Soc. Japan, 66 (1997), 51-59
%
\bibitem{kas} R. Kemmoku and S. Saito,
hep-th 9510007; Jour. Phys. Soc. Japan,
65 (1996), 1881-1884
%
\bibitem{kbe} R. Kemmoku and S. Saito,
Phys. Lett. B, 319 (1993), 471-477; Jour. Phys. A, 29 (1996), 4141-4148
%
\bibitem{kzq} R. Kemmoku and H. Saito,
hep-th 0007122
%
\bibitem{kah} S. Kharchev, A. Mironov, and A. Morozov,
hep-th 9501013
%
\bibitem{kap} B. Khesin, V. Lyubashenko, and C. Roger,
Jour. Fnl. Anal., 143 (1997), 55-97
%
\bibitem{ka} A. Klimyk and K. Schm\"udgen,
Quantum groups and their representations, Springer, 1997
%
\bibitem{kt} B. Kupershmidt,
Lett. Math. Phys., 20 (1990), 19-31
%
\bibitem{kag} B. Kupershmidt and P. Mathieu,
Phys. Lett. B, 227 (1989), 245-250
%
\bibitem{kaz} B. Kupershmidt,
Discrete Lax equations and differential-difference calculus,
Ast\'erisque, 1985
%
\bibitem{kkt} B. Kupershmidt,
Lett. Math. Phys., 20 (1990), 19-31
%
\bibitem{kku} B. Kupershmidt,
Elements of superintegrable systems, Reidel, 1987
%
\bibitem{mzm} J. Madore, S. Schraml, P. Schupp, and J. Wess,
hep-th 0001203
%
\bibitem{man} Y. Manin,
Jour. Sov. Math., 11 (1979), 1-122
%
\bibitem{mzt} Y. Manin,
Quantum groups and noncommutative geometry, CRM Montr\'eal, 1988; Ann. Inst. Fourier,
371 (1987), 191-205
%
\bibitem{mpb} F. Martinez-Moras and E. Ramos,
hep-th 9206040
%
\bibitem{mpa} J. Mas and M. Seco,
Jour. Math. Phys., 37 (1996), 6510-6529
%
\bibitem{mss} J. Matsukidaira, J. Satsuma, and W. Strampp,
Jour. Math. Phys., 31 (1990), 1426-1434
%
\bibitem{mzz} P. vanMoerbeke,
Lectures on integrable systems, World Scientific, 1994, pp. 163-267
%
\bibitem{mal} M. Mulase,
Inter. Jour. Math., 1 (1990), 293-
%
\bibitem{oh} O. Ogievetsky,
Lett. Math. Phys., 24 (1992), 245-255
%
\bibitem{ow} A. Orlov and P. Winternitz,
solv-int 9701008
%
\bibitem{ot} A. Orlov and E. Schulman,
Lett. Math. Phys., 12 (1986), 171-179; Teor. Mat. Fizika, 64 (1985), 323-328
%
\bibitem{smj} R. Sasaki and I. Yamanaka, Comm. Math. Phys., 108 (1987), 691-704;
Advanced Studies Pure Math., 16 (1988), 271-296
%
\bibitem{sz} A. Schwarz,
Mod. Phys. Lett. A, 6 (1991), 611-616; 2713-2725
%
\bibitem{sh} I. Strachan, 
Jour. Phys. A, 28 (1995), 1967-1975; Jour. Geom. Phys., 21 (1997), 255-278;
hep-th 9606101
%
\bibitem{tz} K. Takasaki and T. Takebe,
Inter. Jour. Mod. Phys. A, Supp. 1992, pp. 889-922; Rev. Math. Phys., 7 (1995), 743-808
%
\bibitem{tp} F. Treves,
Introduction to pseudodifferential operators and Fourier integral operators, Vols.
1 and 2, Plenum Press, 1980
%
\bibitem{tu} M. Tu,
Lett. Math. Phys., 49 (1999), 95-103; hep-th 0103083
%
\bibitem{va} A. Volkov,
hep-th 9509024
%
\bibitem{wq} H. Wachter and M. Wohlgenannt,
hep-th 0103120
%
\bibitem{za} M. Zyskin,
hep-th 9504041
%









\end{thebibliography}
\end{document}